\input amstex

\documentstyle{amsppt}
\input xypic
\pageheight{45pc}
\NoBlackBoxes

\define\Q{\Bbb Q}
\define\G{\Bbb G}
\define\Z{\Bbb Z}
\define\N{\Bbb N}
\define\F{\Bbb F}

\define\cal{\Cal}

\magnification 1200
\topmatter
\rightheadtext{Support problem,}
\title  Support problem for the intermediate
jacobians of $l$-adic representations
\endtitle
\author G. Banaszak, W. Gajda, P. Kraso\'n
\endauthor
\address
Department of Mathematics, Adam Mickiewicz University,
Pozna\'{n}, Poland
\endaddress
\email BANASZAK\@math.amu.edu.pl
\endemail
\address
Department of Mathematics, Adam Mickiewicz University,
Pozna\'{n}, Poland
\endaddress
\email GAJDA\@math.amu.edu.pl
\endemail
\address
current: Max Planck Institut f{\" u}r Mathematik in Bonn, Bonn, Germany
\endaddress
\email
GAJDA\@mpim-bonn.mpg.de
\endemail
\address
Department of Mathematics, Szczecin University,
Szczecin, Poland
\endaddress
\email
KRASON\@uoo.univ.szczecin.pl
\endemail

\abstract
We consider the support problem of Erd{\" o}s in the context of
$l$-adic representations of the absolute Galois group of a number field.
Main applications of the results of the paper concern Galois cohomology
of the Tate module of abelian varieties with real and complex
multiplications, 
the algebraic $K$-theory groups of number fields and the integral homology of
the general linear group of rings of integers. We answer
the question of Corrales-Rodrig{\' a}{\~ n}ez and Schoof concerning the
support problem for higher dimensional abelian varieties.
 \endabstract
\endtopmatter

\document

\bigskip

\noindent
\subhead  1. Introduction
\endsubhead

\noindent The support problem for ${\Bbb G}_m$ was first stated by P{\' a}l
Erd{\" o}s who in 1988 raised the following question:
\roster
\item"{}"
{\it let $Supp (m)$ denote the
set of prime divisors of the integer $m.$ Let $x$ and $y$ be two
natural numbers. Are the following two statements equivalent ?
\endroster
\roster
\quad \item"{(1)}"\quad $Supp (x^n {-} 1) = Supp (y^n {-} 1)$\quad for
every\quad $n \in \N,$
\endroster
\roster
\quad \item"{(2)}"\quad $x = y.$}
\endroster

\noindent This question, along with its extension to all number
fields, and also its analogue for elliptic curves, were solved by
Corrales-Rodrig{\' a}{\~ n}ez and Schoof in the paper [C-RS].
Other related support problems can be found in [Ba] and [S]. In
the present paper we investigate the support problem in the
context of $l$-adic representations $$\rho_l\colon\, G_F\to
GL(T_l).$$ The precise description of the class of representations
which are considered is rather technical. It is given by
Assumptions I, II in sections 2 and 3, respectively. This class of representations 
contains powers of the cyclotomic character,
Tate modules of abelian varieties of nondegenarate CM type, and
also Tate modules of some abelian varieties with  real
multiplications (cf. Examples 3.3-3.7).
\medskip

Consider
 the reduction map $$ r_v\colon H^1_{f,S_l}(G_F;\,
T_l)\rightarrow H^1(g_v;\, T_l),$$ for all $v \notin S_l,$ which
is defined on the subgroup $H^1_{f,S_l}(G_F;\, T_l)$ of the Galois
cohomology group $H^1(G_F;\, T_l)$ (see Definition 2.2). 
 $S_l$ denotes here  a finite set of primes which contains primes over $l$
in $F.$ Let $B(F)$ be a finitely generated abelian group such that
for every $l$ there is an injective homomorphism
$$\psi_{F,l}\,\colon\, B(F)\otimes \Z_l \rightarrow
H^1_{f,S_l}(G_F;\, T_l).$$ Let $P$ and $Q$ be two nontorsion
elements of $B(F).$ Put ${\hat P} = \psi_{F,l}(P\otimes 1)$ and
${\hat Q} = \psi_{F,l}(Q\otimes 1).$ Our main point of interest is
the following support problem.
\bigskip

\noindent
\proclaim{Support Problem}\newline
 Let ${\cal P}^{\ast}$ be an infinite set of prime
numbers. Assume that for every $l \in {\cal P}^{\ast}$ the following
condition holds in the group $H^1(g_v;\, T_l):$
\roster
\item"{}" for every integer $m$ and for almost every $v\not\in S_l$
$$
m\,r_v({\hat P})=0\quad\quad \text{implies}\quad\quad
m\,r_v({\hat Q})=0.
$$
\endroster
How are the elements $P$ and $Q$ related in the group $B(F)$ ?
\endproclaim
\bigskip

\subsubhead\rm{\bf  Main results}
\endsubsubhead

\noindent
Let ${\cal P}^{\ast}$ be the infinite set of prime numbers
which we define precisely in section 4 of the paper. We prove the following
theorem.
\bigskip

\proclaim{Theorem A}[Th. 5.1]\newline
Assume that for every
$l \in {\cal P}^{\ast},$
for every integer $m$ and for almost every $v\not\in S_l$
the following condition holds in the group $H^1(g_v;\, T_l):$
\roster
\item"{}"
$$
m\,r_v({\hat P})=0\quad\quad \text{implies}\quad\quad
m\,r_v({\hat Q})=0.$$
\endroster

\noindent Then there exist $a \in {\Bbb Z} -\{0\}$ and $f \in
{\cal O}_E -\{0\}$ such that  $aP + fQ = 0$ in $B(F).$ Here ${\cal
O}_E$ denotes the ring of integers of the number field $E$
associated with the representation $\rho_l$ (see Definition 3.2).
\endproclaim

\noindent
In order to prove Theorem A we investigate 
representations with special properties formulated in Assumption I
and Assumption II. We introduce the notion of the
Mordell-Weil ${\cal O}_E$-module for such representations. Proof of
Theorem A is based on a careful study of reduction maps in Galois
cohomology associated with the given $l$-adic representation satisfying 
 Assumptions I and II. We managed to extend the method of 
[C-RS] to the context of such $l$-adic representations. The main point in the proof
  is to control the relation between  arithmetical properties of the 
images of $l$-adic  representations and reduction maps. Key theorems on the image of 
the representations are proved in the separate work cf. Theorem A and Theorem B of [BGK1]. 
In section 6 of the paper we derive the following 
corollaries of Theorem~A. 
\bigskip

\proclaim{Theorem B}[Cor. 6.4]\newline
Let $P,\, Q$ be two nontorsion elements of
the algebraic $K$-theory group $K_{2n+1}(F),$ where $n$ is an even,
positive integer. Assume that for almost every prime 
\, $v$ 
of \, ${\cal O}_{F}$ and every integer
$m$ the following condition holds in the group  $K_{2n+1}(k_v):$

$$m\,r_v(P)=0 \quad\quad \text{implies}\quad\quad m\,r_v(Q)=0,$$

\noindent
where in this case, $r_v$ is the map induced on the Quillen K-group by the
reduction at $v.$ Then the elements $P$ and $Q$ of $K_{2n+1}(F)$ are linearly dependent
over $\Z.$
\endproclaim
\bigskip

\noindent Note that Theorem B has already been proven by a different method
in [BGK]. 
\bigskip

\noindent
Theorem A has the following corollary concerning the class of
abelian varieties mentioned in the beginning of this Introduction.
\bigskip

\proclaim{Theorem C}[Cor. 6.11]\newline  Let $A$ be an abelian variety
of dimension $g \geq 1$, defined over the number field $F$ and
such that $A$ satisfies one of the following conditions:

\roster
\item"{(1)}" $A$ has the nondegenerate CM type with $End_{F}(A)\otimes \Q$
equal to a CM field $E$ such that $E^{H}\subset F $ (cf. example 3.5)
\item"{(2)}"  $A$ is simple, principally polarised  with  real multiplication by
 a totally real
field $E=End_{F}(A)\otimes \Q$ such that  $E^{H}\subset F $ 
 and the field $F$ is sufficiently large. 
We also assume that $dim\, A = he,$ where
$e=[E:\Q]$ and $h$ is odd (cf. example 3.6) or
$A$ is simple, principally polarised such that
$End_{F}(A) = \Z$ and
$dim\, A$ is equal to $2$\, or \, $6$ (cf. example 3.7 (b)).
\endroster

\noindent
Let $P,\, Q$ be two nontorsion elements of
the group $A(F).$ Assume that for almost every prime 
\, $v$ 
of \, ${\cal O}_{F}$ and
for every integer $m$ the following condition holds in $ A_v(k_v)$
$$m\,r_v(P)=0\quad\quad
\text{implies}\quad\quad m\,r_v(Q)=0.$$

\noindent
Then there exist $a \in \Z -\{0\}$ and $f \in {\cal O}_E -\{0\}$
such that $aP +fQ = 0$ in $A(F).$
\endproclaim
\medskip

\noindent There are two important special cases of
abelian varieties  satisfying conditions of (2) of Theorem C: 
abelian varieties $A$ with $End_{F} (A) = \Z$ such that
$dim\, A $ is an odd integer 
[Se1] (cf. example 3.7 (b)) and
abelian varieties with
real multiplication by a totally real number field
$E = End_{F} (A) \otimes \Q,$  such that $e = g$ 
[R1] (cf. example 3.7 (a)). Note that 
for these abelian varieties the analogues of the 
open image theorem of Serre have been proven  [R1] and [Se1]. The proof of 
Theorem C relies on the analysis 
of the image of the corresponding Galois representation. 
The necessary information on the image of Galois representations on 
$l$- torsion points of abelian varieties from (1) and (2) of Theorem C is 
provided by Theorem 2.1 and Theorem 3.5 of [BGK1]. 
It is worth  mentioning that Theorem C given above provides an answer to the question of
Corrales-Rodrig{\' a}{\~ n}ez and Schoof about the support problem
for higher dimensional abelian varieties [C-RS], p. 277.  
\bigskip

\noindent 
{\it Acknowledgements}:\quad We thank 
the ICTP in Trieste, the MPI in Bonn, the SFB 343 in Bielefeld, 
CRM in Barcelona, Mathematics Departments of the Ohio State University, 
Northwestern University and the Isaac Newton Institute in Cambridge for 
hospitality and financial support during our visits to the 
institutions, while the work on this paper continued. We would like to thank 
the anonymous referee for usefull remarks and suggestions which improved the 
exposition. The research 
was partially financed by a KBN grant 2 PO3A 048 22.
\bigskip

\noindent
\rm {\bf   Notation \rm}

\roster
\item"{$l$}"\quad is an odd prime number.
\item"{$F$}" \quad is a number field, ${\cal O}_{F}$ its ring of integers.
\item"{$G_F$}"\quad $=\quad G({\bar F}/F)$
\item"{$v$}"\quad denotes a finite prime of ${\cal O}_{F}.$
\item"{${\cal O}_{F,S}$}"\quad is the ring of $S-$integers in $F,$
for a finite set
$S$ of prime ideals in ${\cal O}_F$
\item"{$G_{F,S}$}"\quad is the Galois group of 
the  maximal  extension of $F$ 
unramified outside $S$
\item"{$F_v$}"\quad is the completion of $F$ at $v$ and $k_v$ denotes
the residue field ${\cal O}_{F}/v$
\item"{$G_v$}" \quad $ =\quad G({\overline F_v}/F_v)$
\item"{$I_v$}" \quad is the inertia subgroup of
$G_v$
\item"{$g_v$}" \quad $=\quad G({\overline k_v}/k_v)$
\item"{$T_l$}"\quad denotes a free ${\Z_l}$-module of finite rank $d.$
\item"{$V_l$}" \quad $=\quad T_l\otimes_{\Z_l} \Q_l$
\item"{$A_l$}" \quad $= \quad V_l/T_l$
\item"{$\rho_l:$}"\quad $G_F\rightarrow GL(T_l)$ is a representation
unramified outside a fixed finite set $S_l$ of primes of
${\cal O}_{F},$ containing all primes over $l.$
\item"{${\overline{\rho_l}}$}" \quad denotes the residual representation
 $G_F \rightarrow GL(T_l/l)$ induced by $\rho_l.$
\item"{$F_l$}"\quad $=\quad F(A_l[l])$ denotes the number field
${\bar F}^{ker \overline{\rho_l}}.$
\item"{$G_l$}"\quad $=\quad G(F_l/F);$ observe that
$G_l \cong {\overline{\rho_l}}(G_F)$
 is isomorphic to a subgroup of $GL(T/l) \cong GL_d(\Z/l).$
Let $L/F$ be a finite extension and $w$
a finite prime in $L.$ To indicate that $w$ is not over any prime in
$S_l$ we will write $w\notin S_l,$ slightly abusing notation.
\item"{$[H, H]$}" \quad denotes the commutator subgroup of an abstract group
$H.$
\item"{$H^i(G,M)$}" \quad denotes Galois cohomology group of the G-module $M.$
\endroster
\bigskip

\subhead 2. $l$-adic Intermediate Jacobians
\endsubhead

\proclaim{Definition 2.1}  Define
$$H^1_f(G_F; T_l), \quad (\, resp.\,\,\, H^1_f(G_F; V_l)\,)$$
to be the kernel of the natural map:
$$H^1(G_F; T_l) \rightarrow \prod_v H^1(G_v; T_l)/H^1_f(G_v; T_l)$$
$$(\,resp.\,\,\, H^1(G_F; V_l) \rightarrow \prod_v
H^1(G_v; V_l)/H^1_f(G_v; V_l) \,)$$

\noindent
where $H^1_f(G_v; T_l) = i_v^{-1} H^1_f(G_v; V_l)$ via the natural map
 $$i_v : H^1(G_v; T_l) \rightarrow H^1(G_v; V_l).$$
The group $H^1_f(G_v; V_l)$ is defined in [BK] p. 353 (see also [F] p. 115)
as follows:

$$H^1_f(G_v; V_l)  = \cases Ker\, (\, H^1(G_v; V_l) \rightarrow
H^1(I_v; V_l)\,)  & \text{if $v \nmid l$}\\
\\
Ker\, (\, H^1(G_v; V_l) \rightarrow
H^1(G_v; V_l \otimes_{\Q_l} B_{cris})\,) & \text{if $v \mid l$}\, ,
\endcases$$
where $B_{cris}$ is the ring defined by Fontaine (cf. [BK] p. 339).
\endproclaim
\medskip

\noindent
We have the natural maps
$$H^1_f(G_F; T_l) \rightarrow \prod_v H^1_f(G_v; T_l),$$
$$H^1_f(G_F; V_l) \rightarrow \prod_v H^1_f(G_v; V_l).$$
\medskip

\noindent
\proclaim{Definition 2.2} We also define
$$H^1_{f,S_l}(G_F; T_l) \quad (\, resp.\,\,\, H^1_{f,S_l}(G_F; V_l)\, )$$
as the kernel of the natural map:
$$H^1(G_F; T_l) \rightarrow \prod_{v\notin S_l}
H^1(G_v; T_l)/H^1_f(G_v; T_l)$$
$$(\, resp.\,\,\, H^1(G_F; V_l) \rightarrow \prod_{v\notin S_l}
H^1(G_v; T_l)/H^1_f(G_v; V_l)\,).$$ Here $S_l$ denotes a fixed finite set of 
primes of ${\cal O}_F$ containing primes over $l$ and such that the 
representation $\rho_l$ is unramified outside of $S_l.$ 
\endproclaim
\medskip

\noindent
Obviously
$$H^1_{f}(G_F; T_l) \subset H^1_{f,S_l}(G_F; T_l)\quad
\text{and}\quad H^1_{f}(G_F; V_l) \subset H^1_{f,S_l}(G_F; V_l).$$

\noindent Below we define various intermediate
Jacobians associated with the representation $\rho_l,$ (cf. [Sc],
chapter 2).
\medskip

\noindent
\proclaim{Definition 2.3} We put
\roster
\medskip

\item[{1}] $$J(T_l) = \varinjlim_{L/F} H^1(G_L; T_l),
\quad
J(V_l) = \varinjlim_{L/F} H^1(G_L; V_l)$$

\item[{2}] $$J_f(T_l) = \varinjlim_{L/F} H^1_{f}(G_L, T_l), \quad
J_f(V_l) = \varinjlim_{L/F} H^1_{f}(G_L; V_l)$$

\item[{3}] $$J_{f,S_l}(T_l) = \varinjlim_{L/F} H^1_{f,S_l}(G_L; T_l),
\quad
J_{f,S_l}(V_l) = \varinjlim_{L/F} H^1_{f,S_l}(G_L; V_l)$$
\endroster
where the direct limits are taken over all finite extensions $L/F$ of the
number field $F,$ which are  contained in some fixed algebraic closure 
${\bar F}.$
\endproclaim
\medskip

\noindent
\remark{Remark 2.4}
Observe that the groups $J(V_l),$ $J_{f}(V_l)$ and $J_{f,S_l}(V_l)$
are vector spaces over $\Q_l.$
\endremark
\medskip

\noindent \remark{Remark 2.5} Note that we also could have  defined
the intermediate Jacobians of the module $T_l$ for the cohomology
groups of $G_{F,\Sigma}$ for any $\Sigma$ containing $S_l.$
However, if $H^0(g_v; A_l(-1))$ is finite for all $v \notin S_l,$
(as it often happens for interesting examples of $T_l$), then
$$H^1(G_{F,\Sigma}; T_l) = H^1(G_F; T_l).$$
\endremark
\medskip

\noindent
In the sequel we will only consider $l$-adic representations which satisfy 
the following condition. 
\medskip

\proclaim{Assumption I} Assume that for each $l,$ each
finite extension $L/F$ and any prime $w$ of ${\cal O}_L,$ such
that $w\not\in S_l,$ we have $$T_l^{Fr_w} = 0,$$ (or equivalently
$V_l^{Fr_w} = 0$), where $Fr_w\in g_w$ denotes the arithmetic
Frobenius at $w.$
\endproclaim
\medskip

\noindent \remark{Example 2.6} Let $X$ be a smooth projective
variety defined over a number field $F$ with good reduction
at primes $v \notin S_l.$ Let ${\cal X}$ be the
regular, proper model of $X$ over ${\cal O}_{F,S_l}$ and let
${\cal X}_v$ be its reduction at the prime $v$ of ${\cal
O}_{F,S_l}.$ Put ${\overline X} = X \otimes_F{\overline F}$ and
${\overline {\cal X}_v} = {\cal X}_v \otimes_{k_v} {\overline
k_v}.$ In the case when $H^{i}_{et}({\overline X}; \Z_l(j))$ is torsion
free for some $i, j$ such that $i \not= 2j$ we put  $$T_l =
H^{i}_{et}({\overline X}; \Z_l(j)).$$ By the theorem of proper and
smooth base change ([Mi1] VI, Cor. 4.2) there is a natural
isomorphism of $G_v$-modules 
$$H^{i}_{et}({\overline X}; \Z_l(j))
\cong H^{i}_{et}({\overline {\cal X}_v}; \Z_l(j)).\tag{2.7}$$ (cf.
[Ja] p. 322). Since the inertia subgroup $I_v \subset G_v$ acts
trivially on $H^{i}_{et}({\overline {\cal X}_v}; \Z_l(j)),$ we
observe by (2.7) that the representation $\rho_l\,\colon\, G_F
\rightarrow GL(T_l)$ is unramified outside $S_l.$ It follows by
the theorem of Deligne [D1] (proof of the Weil conjectures, see also
[Har] Appendix C, Th. 4.5) that for an ideal $w$  of ${\cal O}_L$
such that $w\not\in S_l,$ the eigenvalues of $Fr_w$ on the vector
space $$V_l = H^{i}_{et}({\overline X}; \Q_l(j))$$ are algebraic
integers of the absolute value $N(w)^{-i/2 + j},$ where $N(w)$
denotes the absolute norm of $w.$ It follows that $T_l^{Fr_w} =
0.$ In the special case when $X=A$ is an abelian variety defined
over $F,$ we have $$T_l = H^{i}_{et}({\overline A}; \Z_l(j)) \cong
\wedge^i H^{1}_{et}({\overline A}; \Z_l) (j) \cong \wedge^i
Hom_{\Z_l}(T_{l}(A) ; \Z_l) (j),$$ which is a free $\Z_l$ module
of rank ${2g \choose i}$ by [Mi2] Th. 15.1. In this paper, most of
the time we will consider the representation $$\rho_l\,\colon\,
G_F \rightarrow GL(T_l (A))$$ of the Galois group $G_F$ on the
Tate module $T_l (A) = H^{1}_{et}({\overline A}; \Z_l)^{\ast}$ of
the abelian variety A defined over $F$, where  for a $\Z_l$-module
$M,$ we put $M^{\ast} = Hom_{\Z_l}(M; \Z_l).$ By the above
discussion we see that $\rho_l$ satisfies Assumption I.
\endremark
\medskip

\proclaim{Lemma 2.8} For every prime $w$ of ${\cal O}_L$ which is not over
primes in $S_l,$ we have:
\roster
\item[{1}] the natural map
$H^1(G_w; T_l)/H^1_f \rightarrowtail  H^1(G_w; V_l)/H^1_f$
is an imbedding,
\item[{2}] $H^1_f(G_w; T_l)_{tor} =  H^1(G_w; T_l)_{tor} = H^0(G_w;
A_l) = H^0(g_w; A_l)$
\item[{3}] $H^1_f(G_w; T_l) =  H^1(g_w; T_l).$
\endroster
\endproclaim
\medskip

\noindent
\demo{Proof}
First part of the lemma is obvious from the definition of $H^1_f(G_w; T_l).$
The second part follows immediately from the first part and
 the diagram (2.9). Note that 
$H^1(G_w; V_l)/H^1_f(G_w; V_l)$ is a $\Q_l$-vector space. To
prove the third part consider consider again the diagram (2.9).

$$\CD
H^0(g_w;\,\,A_l ) @>>> H^1(g_w; \,\,T_l) @>>>
 H^1_f(G_w; \,\,V_l)\\
@V{=}VV @V{}VV @V{}VV\\
H^0(G_w;\,\,A_l )@>>> H^1(G_w;\,\,T_l ) @>>>
H^1(G_w;\,\,V_l )\\
@V{}VV @V{}VV @V{}VV\\
H^0(I_w;\,\,A_l )@>{0}>> H^1(I_w;\,\,T_l ) @>>>
H^1(I_w;\,\,V_l )
\endCD
\tag{2.9}$$
The horizontal rows are exact. The middle and the right vertical columns
are also exact. The left bottom horizontal arrow is zero because $I_w$
acts on $T_l,$ $V_l$ and $A_l$ trivially by assumption. This gives the
exactness of the following short exact sequence.

$$0 \rightarrow H^0(I_w; T_l) \rightarrow H^0(I_w; V_l) \rightarrow
H^0(I_w; A_l) \rightarrow 0$$
In addition  because of Assumption I we have
$$H^0(g_w; V_l) = H^0(G_w; V_l) = 0.$$
Therefore the left upper and middle horizontal arrows are imbeddings.
The right, upper horizontal arrow is defined because of the commutativity
of the lower, right square in the diagram. The middle vertical column
is the inflation restriction sequence. It is actually inverse limit on
coefficients of the inflation-restriction sequence but it remains 
exact with infinite coefficients because we deal with $H^1.$ Now the 
claim follows by diagram chasing.
\qed\enddemo
\medskip

\noindent \remark{Remark 2.10} Observe that the Assumption I
implies, that $H^0(g_w; A_l)$ and $H^1(g_w; T_l)$ are finite for
all $w\not\in S_l.$
\endremark
\medskip

\noindent
\proclaim{Lemma 2.11} For any finite extension $L/F$ the following equalities
hold.

$$H^1_{f, S_l}(G_L; T_l)_{tor} =  H^1(G_L; T_l)_{tor} = H^0(G_L; A_l)$$
\endproclaim
\medskip

\noindent
\demo{Proof} The first equality follows from Lemma 2.8 and the exact sequence.

$$0 \rightarrow H^1_{f, S_l}(G_L; T_l) \rightarrow H^1(G_L; T_l) \rightarrow
\prod_{w \notin S_l} H^1(G_w; T_l)/H^1_f(G_w; T_l).$$
Consider the exact sequence (see [T], p. 261):

$$
\CD
H^0(G_L; V_l) @>>> H^0(G_L; A_l) @>{\partial_L}>> H^1(G_L; T_l).
\endCD
$$

\noindent
By Assumption I we get $H^0(G_L; V_l) = 0.$
Hence by [T], Prop. 2.3, p. 261

$$H^0(G_L; A_l) = H^1(G_L; T_l)_{tor}.$$ Thus, the second equality in
the statement of Lemma 2.11 also holds . \qed\enddemo
\medskip

\noindent
For $w \notin S_l$ consider the following commutative diagram

$$
\CD
H^1_{f, S_l}(G_L; T_l)@>>> H^1(g_w;\, T_l)\\
@AAA @AAA\\
H^0(G_L; A_l)@>>> H^0(g_w;\, A_l).
\endCD
\tag{2.12}$$
The bottom horizontal arrow is obviously an injection.
Hence, by Lemmas 2.8 and 2.11, we obtain the following:
\medskip

\noindent
\proclaim{Lemma 2.13} For any finite extension $L/F$ and  any prime
$w \notin S_l$ in ${\cal O}_{L}$ the natural map

$$
\CD
r_{w}\,\,:\,\,H^1_{f, S_l}(G_L; T_l)_{tor}@>>> H^1(g_w;\, T_l)
\endCD
$$
is an imbedding.
\endproclaim
\medskip

\noindent
\proclaim{Proposition 2.14} We have the following exact sequences

$$0 \rightarrow A_l \rightarrow J(T_l) \rightarrow J(V_l)
\rightarrow 0$$
$$0 \rightarrow A_l \rightarrow J_{f,S_l}(T_l) \rightarrow J_{f,S_l}(V_l)
\rightarrow 0$$
In particular
$$J(T_l)_{tor} = J_{f,S_l}(T_l)_{tor} = A_l$$
and the groups
$$J(T_l) \quad and \quad J_{f,S_l}(T_l)$$
are divisible.
\endproclaim
\medskip

\noindent
\demo{Proof} Consider the following long exact sequence (see [T] p. 261)

$$H^0(G_L;\,A_l) \rightarrow H^1(G_L;\,T_l) \rightarrow
H^1(G_L;\,V_l) \rightarrow H^1(G_L;\,A_l).$$
Taking direct limits with respect to finite extensions $L/F$
gives the following short exact sequence.
$$0 \rightarrow A_l \rightarrow J(T_l) \rightarrow
J(V_l) \rightarrow 0$$
This short exact sequence fits into the following commutative diagram

$$\eightpoint
\CD
0 @>>>{}0 @>>> \varinjlim_{L/F} \prod_{w \notin S_l} H^1(G_w; \,\,T_l)/H^1_f @>>>
\varinjlim_{L/F} \prod_{w \notin S_l} H^1(G_w; \,\,V_l)/H^1_f\\
@. @A{}AA @A{}AA @AAA\\
0@>>> A_l@>>> J(T_l)@>>> J(V_l)@>>> 0\\
@. @A{=}AA @A{}AA @AAA\\
0@>>> A_l@>>> J_{f,S_l}(T_l)@>>>  J_{f,S_l}(V_l)@>>> 0\\
@. @AAA @AAA @AAA\\
@. 0 @. 0 @. 0\\
\endCD
\tag{2.15}$$
The rows and columns of the diagram are exact. The exactness on the right of
the bottom horizontal sequence follows from the injectivity of the top,
nontrivial, horizontal arrow by Lemma 2.8.
\qed\enddemo
\medskip

\noindent
\proclaim{Proposition 2.16} Let $L$ be a finite extension
of $F.$ Then we have isomorphisms:

\roster
\item[{1}] $H^1(G_L; T_l) \cong J(T_l)^{G_L},$
\item[{2}] $H^{1}_{f,S_l}(G_L; T_l) \cong J_{f,S_l}(T_l)^{G_L}.$
\endroster
\endproclaim

\noindent
\demo{Proof} Under condition of Assumption I
the proof of claim (1) is done in the
same way as the proof of (4.1.1) of [BE]. To prove (2) take an arbitrary
finite Galois extension $L^{\prime}/L$ and consider the
following commutative diagram.

$$
\CD
 0 @. 0\\
@VVV @VVV\\
H^{1}_{f,S_l}(G_L; T_l) @>>>
H^{1}_{f,S_l}(G_{L^{\prime}}; T_l)^{G(L^{\prime}/L)}\\
@VVV @VVV\\
H^{1}(G_L; T_l) @>>{\cong}>
H^{1}(G_{L^{\prime}}; T_l)^{G(L^{\prime}/L)}\\
@VVV @VVV\\
\prod_{w \notin S_l} H^1(I_w; \,\,T_l) @>>>
\prod_{w \notin S_l}\bigl
(\prod_{w^{\prime}|w} H^1(I_{w^{\prime}}; \,\,T_l)\bigl)^{G(L^{\prime}/L)}
\endCD
\tag{2.17}$$ The columns of this diagram are exact. The upper
horizontal arrow is trivially an imbedding. 
The middle horizontal arrow is an isomorphism. This follows directly from
claim (1). Since the representation $\rho_l$ is unramified outside 
$S_l$ then using Th. 8.1 and Cor. 8.3 Chap. I of [CF] and Kummer pairing
we get the following commutative diagram

$$
\CD
H^1(I_w; \,\,T_l) @>>> H^1(I_{w^{\prime}}; \,\,T_l)\\
@VV=V @VV=V\\
Hom_{cts} (I_w; \,\,T_l) @>>>
Hom_{cts} (I_{w^{\prime}}; \,\,T_l)\\
@VV{\cong}V @VV{\cong}V\\
Hom_{cts} (\Z_l(1); \,\,T_l)@>>>
Hom_{cts} (\Z_l(1); \,\,T_l)
\endCD
\tag{2.18}$$
Since $L_{w^{\prime}}/L_{w}$ is a finite extension, the bottom horizontal 
arrow is induced by a nontrivial (hence injective) 
homomorphism of $\Z_l$-modules $\Z_l(1)\rightarrow \Z_l(1).$
Because $T_l$ is a free $\Z_l$-module, every nontrivial
homomorphism of $\Z_l$-modules $\Z_l(1)\rightarrow T_l$ 
is injective. Hence the bottom horizontal arrow in the diagram 
(2.18) is injective.
So the bottom horizontal arrow in diagram (2.17) is also an imbedding.
Now claim (2) follows by taking  direct limits over
$L^{\prime}$ in diagram (2.17) and chasing the resulting diagram. \qed\enddemo

At the end of this section we give some additional information
about the reduction map $$ r_v\colon H^1_{f,S_l}(G_F;\,
T_l)\rightarrow H^1(g_v;\, T_l). $$

\proclaim{Proposition 2.19}
Let ${\hat P}\in H^1_{f,S_l}(G_F;\,T_l)$
be a nontorsion element. Given $M_1=l^{m_1}$ a fixed power of $l,$
there exist infinitely many primes $v\not\in S_l$ such that
$r_v({\hat P})\in H^1(g_v;\, T_l)$ is an element of order at
least $M_1.$
\endproclaim

\noindent {\sl Proof.} Let $M$ be a power of $l$ which we will
specify below. Let $F_M$ denote the extension $F(A[M]).$ Consider
the following commutative diagram.

$$ \CD H^{1}_{f,S_l}(G_F;\,T_{l})/M @>{r_v}>>
H^{1}(g_v;\,T_{l})/M  \\ @V{h_1}VV    @VVV \\
H^{1}_{f,S_l}(G_F;\,A[M])  @>{r_v}>>  H^{1}(g_v;\,A[M])\\
@V{h_2}VV    @VVV \\ H^{1}_{f,S_l}(G_{F_M};\, A[M])  @>{r_w}>>
H^{1}(g_{w};\, A[M])\\ @V{h_3}VV    @VVV \\ Hom(G_{F_M};\, A[M])
@>{r_w}>>  Hom(G_w;\, A[M])\\ @V{h_4}V{\cong}V   @VV{=}V
\\ Hom(G_{F_M}^{ab};\, A[M])  @>{r_w}>> Hom(g_w;\, A[M])\\
\endCD
\tag 2.20
$$
The horizontal arrows in the diagram (2.20) are induced by the reduction
maps. We describe the vertical   map.  The map 
$h_2$ is an injection (cf. Prop. 4.3 (1)). The map $h_3$ is the injection which comes
from the long exact sequence in cohomology associated to the
following exact sequence of $G_{F_M}$-modules: $$ \CD 0@>>> A[l]
@>>> J_{f, S_l}(T_l) @>{\times l}>> J_{f, S_l}(T_l)@>>> 0.
\endCD\tag 2.21
$$  The vertical maps on the right hand side of the diagram (2.20)
are defined in the similar way. Consider the nontorsion element
${\hat P}\in H^1_{f,S_l}(G_F;\,T_l).$ Let $l^s$ be the largest
power of $l$ such that ${\hat P}=l^s{\hat R}$ for an ${\hat R}\in
H^1_{f,S_l}(G_F;\,T_l).$ Such an $l^s$ exists since
$H^1_{f,S_l}(G_F;\,T_l)$ is a finitely generated $\Z_l$-module. We
put $M=M_1l^{s}.$ Let ${P^{\prime}}$ be the image of ${\hat P}$ in
$Hom(G_{F_M}^{ab};\, A[M])$ under the composition of the maps
$h_1,$ $h_2,$ $h_3$ and $h_4$. Since the maps $h_1,$ $h_2,$ $h_3$
and $h_4$ are injective, the element ${P^{\prime}}$ is of order
$M_1.$ By the Chebotarev density theorem there exist infinitely
many primes $w\not\in S_l$ such that the map $r_w$ preserves the
order of $P^{\prime}.$ Hence, for those $w$ the element ${\hat P}$
is mapped by the composition of left vertical and lower horizontal
arrows onto an element whose order is $M_1.$ The commutativity of
(2.20) implies that
 $r_v({\hat P})\in H^1(g_v;\,T_l)$ is of order at least $M_1$ for
the primes $v=w\cap {\cal O}_{F,S_l}.$
\qed
\bigskip

\proclaim{Corollary 2.22}\newline \noindent Let ${\hat P}\in
H^1_{f,S_l}(G_F;\, T_l)$ be an element which maps onto a generator
of the free $\Z_l$ - module  $H^1_{f,S_l}(G_F;\, T_l)/tor$. There
exist infinitely many primes $v\not\in S_l$ such that $r_v({\hat
P})$ is a generator of a cyclic summand in the $l$-primary
decomposition of the group $H^1(g_v;\,T_l).$
\endproclaim
\bigskip\bigskip

\subsubhead\rm{\bf 3. Specification of $l$-adic representations }
\endsubsubhead
\bigskip

In addition to Assumption I the representations which we consider are 
supposed to satisfy Assumption II stated below. 
In order to formulate the assumption we introduce more notation. We fix
a finite extension $E/\Q$ of degree $e = [E:\Q]$ such that the
Hilbert class field $E^H$ of $E$ is contained in $F.$ We assume
that  each  prime $l$ splits completely in $F.$ Let
$$(l)=\lambda_1 \dots \lambda_e$$ be the decomposition of the
ideal $(l)$ in ${\cal O}_E.$ We also assume that ${\cal O}_E$ acts on $T_l$
in such a way that $T_l$ is a free ${\cal O}_{E,l} = {\cal
O}_E\otimes \Z_l$
 module of rank $h$ and that the action of ${\cal O}_{E,l}$
commutes with the action of $G_F$ given by the representation
$\rho_l.$ It is clear that $e$ divides $d=dim\, \rho_l$ and $h =
{d\over e}.$ Put $E_l = {\cal O}_{E,l}\otimes_{\Z_l} \Q_l.$ In
addition, we denote by $E_{\lambda_i}$ the completion of $E$ at
$\lambda_i$ and by ${\cal O}_{\lambda_i}$ the ring of integers in
$E_{\lambda_i}.$ Now, it is obvious that $${\cal O}_{E,l} =
\prod_{i=1}^{e} {\cal O}_{\lambda_i}\quad\quad
\text{and}\quad\quad E_l = \prod_{i=1}^{e} E_{\lambda_i}.$$ 
Since
$T_l$ has the ${\cal O}_{E,l}$-module structure, we can represent
$V_l$ and $A_l$ as follows: $$V_l=T_l \otimes_{{\cal O}_{E,l}}
E_l$$ $$A_l=T_l \otimes_{{\cal O}_{E,l}} E_l/{\cal O}_{E,l} =
\oplus_{i=1}^{e} T_l \otimes_{{\cal O}_{E,l}}E_{\lambda_i}/{\cal
O}_{\lambda_i}
 = \oplus _{i=1}^{e} A_{\lambda_i},$$
where we put $A_{\lambda_i}$ $=$
$T_l \otimes_{{\cal O}_{E,l}} E_{\lambda_i}/{\cal O}_{\lambda_i}.$

 Note that every prime ideal $\lambda_i$ is
principal, because by assumption $E^H \subset F.$ Hence,
$\lambda_i = (\pi_i)$ for some $\pi_i \in {\cal O}_E.$ In this
case $E_{\lambda_i}/{\cal O}_{\lambda_i} \cong \Q_l/\Z_l$ for each
$i,$ hence all $A_{\lambda_{i}}$ are divisible groups of the same
corank $h.$ Observe that $A_{l}[{\lambda_i}^k] =
A_{\lambda_i}[{\lambda_i}^k]$ and $$A_l[l]\,\, \cong\,\,
\bigoplus_{i=1}^e\, A_{l}[{\lambda_i}],\tag{3.1}$$ where
$dim_{\Z/l}\, A_l[\lambda_i] = h,$ for all $1 \leq i \leq e.$ By
assumptions and decomposition (3.1) it is clear that the image
$\overline{\rho_{l}}(G_F)$ of the representation
$\overline{\rho_{l}}$ is contained in the subgroup of $GL_d(\Z/l)$
which consists of matrices  of the form $$ \pmatrix
C_1&0&\ldots&0\\
               0&C_2&\ldots&0\\
               \vdots&\vdots&\vdots&\vdots\\
               0&0&\ldots&C_{e}
       \endpmatrix ,
$$
where $C_i \in GL_h(\Z/l),$ for all $1 \leq i \leq e.$
Hence, we can consider the image of $\overline{\rho_{l}}$ as a
subgroup of the product $\prod_{i=1}^{e} GL_h(\Z/l).$
\bigskip

\proclaim{Assumption II}{\rm
Let ${\cal P} = {\cal P}(\rho)$ be an infinite set of
prime numbers $l > 3,$ which split completely in $F$ and such that the 
$l$-adic representation $\rho_l$ satisfies the following
conditions.
\roster
\item"{(1)}" If $h > 1,$ then for each $1 \leq i \leq {e},$ there is
a subgroup $H_i \,\trianglelefteq\, G(F(A[\lambda_i])/F) 
\subset GL_h(\Z/l)\cong GL(A_l[\lambda_i])$
such that:
\roster
\quad \item"{$(i)$}"   the subgroup
$H_1\times\dots \times H_i\times\dots\times
H_e$ of the group $\prod_{i=1}^{e} GL_h(\Z/l)$ is contained in
$Im\, \overline{\rho_{l}} = G_l$ and
$H_1\times\dots \times H_i\times\dots\times
H_e$ has index prime to $l$ in $G_l$
\endroster
\roster
\quad \item"{$(ii)$}" $H_i$ acts irreducibly on
$A_l[\lambda_i] \cong (\Z/l)^h,$
\endroster
\roster
\quad \item"{$(iii)$}" $H_i/[H_i, H_i]$ has order prime to $l,$
\endroster
\roster
\quad \item"{$(iv)$}" there exist matrices $\sigma_i,\, \beta_i \in H_i$
such that $1$ is an eigenvalue of $\sigma_i$ with eigenspace of dimension $1$
and $1$ is not an eigenvalue of $\beta_i,$
\endroster
\roster
\quad \item"{$(v)$}" the centralizer of $H_i$ in $GL_h(\Z/l)$ is
$(\Z/l)^{\times} I_h$ i.e. if $\sigma \in GL_h(\Z/l)$ and
$\sigma\, \gamma = \gamma\, \sigma$ for all $\gamma \in H_i,$ then
$\sigma$ is a scalar matrix,
\endroster
\roster
\quad\item"{$(vi)$}"  for each $1 \leq i \leq e$ the group
$H_i$ contains
a nontrivial subgroup $D_i^0$ of the group
$\{aI_h;\,\, a \in (\Z/l)^{\times}\}
\subset GL_h(\Z/l)$ of scalar matrices.
\endroster
\endroster
\roster
\item"{(2)}" If $h = 1,$ we require that $G_l = G(F_l/F)$ satisfies two
 conditions:
\roster
\quad \item"{$(i)$}" for every $1 \leq i \leq d,$ there is a diagonal matrix
$\sigma_i = diag\, (\mu_1,\dots, \mu_d)$ in the group $G_l$ with
$\mu_i = 1$ and $\mu_j \not= 1,$ for all $j \not= i,$
\endroster
\roster
\quad \item"{$(ii)$}" there is an isomorphism of rings
$\Z/l[G_l] \cong {\cal O}_E/l,$
where $\Z/l[G_l]$ denotes a subring of ${\cal O}_E/l$
generated by $\Z/l$ and the image of $G_l$ in ${\cal O}_E/l$
via the natural imbedding $G_l \rightarrow ({\cal O}_E/l)^{\times}.$
\endroster
\endroster}
\endproclaim

\proclaim{Definition 3.2}{\rm Let $\{B(L)\}_L$ be a direct
system of ${\cal O}_E$-modules indexed by all finite field
extensions $L/F.$ The structure maps of the system are induced by
inclusions of fields. We assume that for every embedding of fields
$L \rightarrow L^{\prime}$ the structure map $B(L) \rightarrow
B(L^{\prime})$ is a homomorphism of ${\cal O}_{E}$-modules. Let us
put $B({\overline F}) = \varinjlim_{L/F}\, B(L).$ Let $\rho$ and
${\cal P}$ be as in Assumption II. The system $\{B(L)\}_{L}$ is
called the Mordell-Weil ${\cal O}_E$-module of the pair $(\rho,\, {\cal P})$ if
the following conditions are satisfied: \roster
\item"{(1)}" $B(L)$ is a finitely generated ${\cal O}_E$-module
for all $L.$
\item"{(2)}" There are functorial homomorphisms of ${\cal O}_E$-modules 
$$\psi_{L,l}\, :\, B(L) \longrightarrow
H^1_{f,S_{l}}(G_L; T_l),$$
where $L$ is as above and $l\in {\cal P},$ such that:
\roster
\quad\item"{$(i)$}" for every $l\in {\cal P},$ the induced map
$$\psi_{L,l} \otimes \Z_l: B(L)\otimes \Z_l
\rightarrow H^1_{f,S_{l}}(G_L; T_l)$$
is an isomorphism or
\endroster
\roster
\quad\item"{$(ii)$}" for every $l\in {\cal P},$ the map
$\psi_{L,l} \otimes \Z_l$ is an imbedding,
the group $B({\overline F})$ is a discrete $G_F$-module which is
divisible by $l,$ and for every $L$ we have:
$B({\overline F})^{G_L} \cong B(L)$ and $H^0(G_L; A_l) \subset B(L).$
\endroster
\endroster
}
\endproclaim

\medskip
We end this section with the examples of Mordell-Weil
${\cal O}_E$-modules related to $l$-adic representations which satisfy
Assumptions I and II.

\medskip
\remark{Example 3.3}{\rm Consider the $l$-adic representation
$$\rho_{l}\,\,:\,\, G({\bar F}/F) \rightarrow GL(\Z_l(1)) \cong
GL_1(\Z_{l}) \cong \Z_l^{\times}$$ given by the cyclotomic
character. In this case $T_l=\Z_l(1),$ $V_l=\Q_l(1)$ and
$A_l=\Q_l/\Z_l(1).$ This representation is given by the action of
$G_F$ on the Tate module of the multiplicative group scheme ${\Bbb
G}_m/F.$  Let $S$ be any finite set of primes in ${\cal O}_F.$
Denote by $S_l$ the set of primes consisting of primes in $S$ and
primes in $F$ over $l.$ Put $B(L) = {\G}_m({\cal O}_{L, S}) =
{\cal O}_{L, S}^{\times}$ for any finite extension $L/F.$ The
Kummer map (which is obviously injective) $$B(L)\otimes \Z_l
\rightarrow H^1(G_{L,S_{l}}; \Z_l(1)) \rightarrow H^1(G_L;
\Z_l(1))$$ factors naturally through $$\psi_{L,l}\,\,:\,\,
B(L)\otimes \Z_l \rightarrow H^1_{f, S_{l}}(G_F; \Z_l(1))$$ In
this case we take $E = \Q$ hence ${\cal O}_E = \Z.$ We take ${\cal
P}$ to be the set of all prime numbers $l$ such that
$G(F(\mu_l)/F)$ is nontrivial.}
\endremark
\medskip

\remark{Example 3.4}{\rm
Let $n$ be a positive integer. Let $T_l=\Z_l(n+1),$ hence
$V_l=\Q_l(n+1)$ and $A_l=\Q_l/\Z_l(n+1).$  Consider the one dimensional
representation $$
\rho_l\,\,:\,\, G_F\rightarrow GL(T_l)\cong \Z_l^{\times}$$
which is given by the $(n+1)$-th tensor power of the cyclotomic
character. For each odd prime number $l$ and for a finite extension 
$L/F$ consider the Dwyer-Friedlander map [DF]
$$ K_{2n+1}(L)\rightarrow
K_{2n+1}(L) \otimes \Bbb Z_l\rightarrow H^1(G_L;\, \Bbb Z_l(n+1)).$$
Let $C_L$ be the subgroup
of $K_{2n+1}(L)$ which is generated by the $l$-parts of kernels of
Dwyer-Friedlander maps for all odd primes $l.$
We define the group $B(L)$ by putting $$B(L)=K_{2n+1}(L)/C_L.$$

\noindent Note that the group $C_L$ is finite by [DF] and it
should vanish if the Quillen-Lichtenbaum conjecture holds.
Note that in this case
$$H^1(G_L;\, \Z_l(n+1))\cong H^1(G_{L,S_l};\, \Bbb Z_l(n+1))  \cong 
H^1_{f,S_l}(G_L;\, \Bbb Z_l(n+1)).$$
It follows by the definition of $B(L)$ and surjectivity of the Dwyer-Friedlander map
 that
$$
\psi_{L,l}\quad\colon\quad B(L)\otimes\Z_l\cong H^1(G_L;\, \Bbb Z_l(n+1)).
$$}
\endremark

In the following three examples we discuss representations which
come from Tate modules of abelian varieties. 
Let $A/F$ be a simple abelian variety
of dimension $d$ over a number field $F.$  As usual, we denote by  $T_l
= T_l(A)$ the Tate module of $A.$ Consider the $l$-adic
representation $$\rho_l\,\,:\,\, G_F \rightarrow GL(T_l(A)).$$
Assumption I holds due to the Weil conjectures (cf. [Sil], pp.
132-134).  Let $S$ be the set of prime ideals of $F$  at which $A$
has bad reduction. By the Kummer pairing and Serre-Tate theorem
([ST], Th. 1, p. 493 and Corollaries 1 and 2 of Manin's Appendix
II to the book [M]) we have a natural imbedding
$$\psi_{L,l}\quad\colon\quad A(L) \otimes \Z_l \rightarrow H^1_{f,
S_l}(G_L; T_l(A)).$$ Put $B(L) = A(L)$ for any finite extension
$L/F.$
\medskip

\remark{Example 3.5} \rm Let  $A/F$ be a simple abelian variety with
complex multiplication by a CM field $E$ (cf. [La]) such that $E^H
\subset F,$ where $E^H$ is the Hilbert class field of $E.$ We
assume that CM type of $A$ is nondegenerate (cf. Def.2.1, [BGK1])
 and defined over $F$ . Condition ($i$) of Assumption II (2)
holds by Theorem 2.1, [BGK1] (for CM elliptic curves it
also follows by an alternative argument cf. [C-RS], Lemma 5.1, p.
286). Condition ($ii$) of Assumption II (2) follows by
Proposition, p. 72 of [R2].  We take ${\cal P}$ to be the set of
prime numbers $l$ which split completely in $F$ and such that $A$
has a good reduction at $l.$
\endremark
\bigskip

\remark{Example 3.6} \rm Consider a simple, principally polarised abelian
variety $A/F$ such that $E = End_{F}(A)\otimes {\Q} = End_{\bar
F}(A)\otimes {\Q}$ (cf. [R1] and [C])  where $e = [E:\, \Q]$ and
$2h\, e = 2g$ with $h$ and odd integer. In addition, we choose $F$
to be a number field satisfying conditions indicated in the discussion 
which follows Theorem 3.1 of [BGK1] and such
that $E^H \subset F.$ We take ${\cal P}$ to be the set of prime
numbers $l \gg 0$ which split completely in $F,$ and such that 
$A$ has  good reduction at $l.$ 
Hence by Theorem 3.5 of [BGK1] we get $$\prod_{i=1}^e
Sp_{2h}(\F_l) =  [G_l,\, G_l].$$

\noindent Taking $H_i = Sp_{2h}(\F_l),$ for all $1 \leq i \leq
e,$ we observe that conditions of Assumption II (1) are fulfilled
since
\roster
\quad \item"{$(i)$}" $\prod_{i=1}^e Sp_{2h}(\F_l)
\subset G_l,$ and the quotient group $GSp_{2h}(\F_l)/
Sp_{2h}(\F_l)$ has order prime to $l,$
\endroster
\roster \quad \item"{$(ii)$}" $Sp_{2h}(\F_l)$ acts   on
$A_l[\lambda_i] \cong (\Z/l)^{2h},$ in an irreducible way.
\endroster
\roster
\quad \item"{$(iii)$}" $Sp_{2h}(\F_l)$ modulo its center is a simple group.
\endroster
\roster
\quad \item"{$(iv)$}"  matrix $\sigma_i \in Sp_{2h}(\F_l)$
$$\sigma_i = \pmatrix J_{h}(1)& J_{h}(1)\\
                      O     &(J_{h}(1)^t)^{-1}
\endpmatrix
$$
has eigenvalue $1$ with the eigenspace of dimension $1$
where
$J_{h}(1)$ is the ${h}\times {h}$ Jordan block matrix with 1 as the eigenvalue
\noindent
and $\beta_i = -I_{2h} \in Sp_{2h}(\F_l)$ does not have 1 as an eigenvalue.
\endroster
\roster
\quad \item"{$(v)$}" The centralizer of $Sp_{2h}(\F_l)$ in
$GL_{2h}(\F_l)$ is
$(\F_l)^{\times} I_{2h}.$
\endroster
\medskip

\noindent Observe that condition (1) (vi) of Assumption II is satisfied
since obviously $-I_{2h} \in Sp_{2h}(\F_l).$
\endremark
\bigskip

\noindent There are two special cases of Example 3.6 that have been
considered extensively in the past.
\medskip

\noindent
\remark{Example 3.7} \rm (a)\,\, Let $A/F$ be a simple, principally polarised abelian variety
with real multiplication by a totally
real field $E = End_{F}(A)\otimes {\Q} = End_{\bar F}(A)\otimes {\Q}$ 
such that $e = g$ and $h = 1$ (cf. [R1]).
We choose $F$ to be such a number field that
$E^H \subset F.$ We take ${\cal P}$ to be the set of prime numbers
$l$ which split completely in $F$ and such that $A$ has a good
reduction at $l.$ Theorem 5.5.2, p. 801, [R1] or Theorem 3.5 of [BGK1]
implies that
the image of the representation ${\bar \rho_l}$ contains the
subgroup
$$\prod_{i=1}^g SL_{2}(\F_l) = \prod_{i=1}^g Sp_{2}(\F_l),$$
therefore the
representation ${\bar \rho_l}$ satisfies Assumption II (1).
\bigskip

(b)\,\, Let $A/F$ be a simple, principally polarised abelian variety with the property
that $End_{\bar F}(A) = \Z$ and $g = dim\, A$ is odd or equal to
$2$ or $6.$ In this case $E = \Q$ hence $e = 1$ and $h = g.$ By the
theorem of Serre ([Se1]  Th. 3) the image of the representation
${\bar\rho_l}$ equals $GSp_{2g}(\F_l)$ (hence contains
$Sp_{2g}(\F_l)$) for almost all $l.$  We take ${\cal P}$ to be the
set of prime numbers such that the image of ${\bar\rho_l}$ equals
$GSp_{2g}(\F_l)$ and $A$ has  good reduction at $l.$ Hence the image
of the representation ${\bar \rho_l}$ satisfies  condition (1)
 of Assumption II.
\endremark
\medskip

\noindent It is rather hard to find further examples of
Mordell-Weil ${\cal O}_E$-modules satisfying condition (2)(i) of
Definition 3.2. Indeed, if we concentrate on finding a Mordell-Weil
${\cal O}_E$-module associated to $T_l$ coming from {\'e}tale
cohomology of a smooth proper scheme $X$ over $F,$ then we should
first prove Conjecture 5.3 (ii) p. 370 of [BK], for such an $X.$
\bigskip\bigskip

\subhead{4. Key Propositions}
\endsubhead

\noindent
\proclaim{Definition 4.1}
Let 
$$\phi_P\,\,\colon \,\, G_{F_l} \rightarrow A_l[l]$$
be the map:

$$\phi_P(\sigma) = \sigma({1 \over l}{\hat P}) - {1 \over l}{\hat P}$$
where $P \in B(F)$ and ${\hat P}$ is the image of $P$ via the natural map
$$B(F)\rightarrow B(F)\otimes \Z_l \rightarrow  H^1_{f,S_{l}}(G_F; T_l)
\subset J_{f, S_l}(T_l)$$
\endproclaim
\medskip

\noindent \remark{Remark 4.2} Note that \, ${1 \over l}{\hat P}$ \,
makes sense in $J_{f, S_l}(T_l)$ since the last group is divisible
due to Proposition 2.14. The element ${1 \over
l}{\hat P}$ is defined up to an element of the group $A_l[l].$
\endremark
\medskip

\noindent
\proclaim{Proposition 4.3} Suppose that the  Assumptions I and II are fulfilled.
Then the following properties hold.
\roster
\item[{1}] $H^r(G(F_l/F); A_l[l])=0$ for $r\geq 0$ and all $l \in {\cal P},$
except the case of trivial $G_l$-module $A_l[l]$ when $r = 0$ and $d = 1.$
\item[{2}] The map $H^1_{f,S_{l}}(G_F; T_l)/l \longrightarrow
H^1_{f,S_{l}}(G_{F_l}; T_l)/l$
is injective for all $l \in {\cal P}.$
\item[{3}] The map $B(F)/lB(F)\longrightarrow B(F_l)/lB(F_l)$
is injective for all $l \in {\cal P}.$
\item[{4}] Let $P\in B(F).$ If $l \in {\cal P}$ does not divide $\sharp\, B(F)_{tor}$
and $P \notin \lambda_{i}B(F)$ for all $1 \leq i \leq e,$ then the
map $\phi_P$ is surjective.
\endroster
\endproclaim
\medskip

\noindent
\demo{Proof}
(1) First let us consider the case $h > 1.$
The group $D^0 = \prod_{j=1}^e D_i^0$ can be regarded as a subgroup of
$G_l$ once we identify $G_l$ with its
image via $\overline{\rho_{l}}.$ $D^0$ is a normal
subgroup of $G_l.$
Assumption II (1) (vi) allows us to consider
the Hochschild-Serre spectral sequence
$$
E_2^{r,s} = H^r(G_l/D^0;\, H^s(D^0;\, A_l[l]))\Rightarrow H^{r+s}(G_l;\,
A_l[l]).\tag{4.4}
$$
Observe that $H^0(D^0;\, A_l[l]) =
\oplus_{i=1}^e H^0(D_i^0;\, A_l[\lambda_i]) = 0$ because by definition
$D_i^0$ is nontrivial and
acts by matrix multiplication (actually scalar multiplication)
on the $\Z/l$ vector space $A_l[\lambda_i] \cong (\Z/l)^h.$
The groups $H^s(D^0;\, A_l[l])$ vanish for
$s >  0,$ since $l$ is odd by assumption and  the order
of $D^0$ is prime to $l.$ Hence the claim (1) follows
for $h > 1.$\,
Now let $h = 1.$
Note that $G_l$ is isomorphic to a subgroup of diagonal matrices in
$GL(A_l[l]) = GL_d(\Z/l)$.
Since $G_l$ has order relatively prime to $l,$
$H^s(G_l;\, A_l[l]) = 0$ for $s > 0.$
It follows easily by Assumption II (2) (i) that
$H^0(G_l;\, A_l[l])  = 0,$
for all $l \in {\cal P}$ and $d > 1.$  This proves (1) in the case $h = 1.$
If $d = 1,$ then $H^0(G_l;\, A_l[l])  = 0$\,\,\, ($ = A_l[l]$ resp.)\,\,
if $A_l[l]$ is nontrivial (trivial resp.)
$G_l$-module.
\medskip

(2) By Prop. 2.14, we have the following
short exact sequence:
$$
\CD
0@>>> A_l[l] @>>> J_{f, S_l}(T_l) @>{l}>> J_{f, S_l}(T_l)@>>> 0.
\endCD
$$

\noindent
By the long exact sequence in cohomology associated to this exact
sequence and Proposition 2.16, we obtain the
commutative diagram in which the horizontal
maps are injections.

$$
\CD
 0 @. 0\\
@VVV @VVV\\
ker\, \alpha @>>> H^1(G_l;\, A_l[l])\\
@VVV @VVV\\
H^1_{f,S_l}(G_F;\, T_l)/l @>>> H^1(G_F;\, A_l[l])\\
@VV{\alpha}V @VV{\gamma}V\\
H^1_{f,S_l}(G_{F_l};\, T_l)/l @>>> H^1(G_{F_l};\, A_l[l])
\endCD
\tag{4.5}$$
However, $ker\, \alpha =0,$ since it injects into the
group $H^1(G_l;\, A_l[l])$ which vanishes by part (1) of the proposition.

(3)  Let us first consider the case (2) (i) of Definition 3.2.
Because the map

$$
B(L)\otimes \Z_l \longrightarrow H^1_{f,S_l}(G_L;\, T_l),
$$

\noindent
is an isomorphism, the group $B(L)/l$ is isomorphic to
$H^1_{f,S_l}(G_L;\, T_l)/l.$
This shows that the horizontal maps in the commutative diagram

$$
\CD
B(F)/l@>>> H^1_{f,S_l}(G_F;\, T_l)/l\\
@VVV @VV{\alpha}V\\
B(F_l)/l@>>> H^1_{f,S_l}(G_{F_l};\, T_l)/l
\endCD
\tag{4.6}$$

\noindent 
are isomorphisms. Since we have proved in (2) that the
map $\alpha$ is an injection, diagram (4.6) gives the claim (3).
Now consider the case (2) ($ii$) of Definition 3.2. We get the
exact sequence of $G_F$-modules:

$$
\CD
0 @>>> A_l[l] @>>> B({\overline F})@>{l}>> B({\overline F}) @>>> 0
\endCD
$$ This gives the following commutative diagram with injective
horizontal arrows:

$$
\CD
 0 @. 0\\
@VVV @VVV\\
ker\, \beta @>>> H^1(G_l;\, A_l[l])\\
@VVV @VVV\\
B(F)/l@>>> H^1(G_F;\, A_l[l]) \\
@VV{\beta}V @VV{\gamma}V\\
B(F_l)/l@>>> H^1(G_{F_l};\, A_l[l])
\endCD
\tag{4.7}$$
Since  by (1) the map $\gamma$ is injective for all $l \in {\cal P}$,
 the map $\beta$ is also injective for all $l \in {\cal P}.$
\medskip

(4) We easily check that the image of the map $\phi_P$ is
$G_F$-invariant. If $\phi_P$ were not surjective, then $Im\,
\phi_P$ would be a proper $G_F$ submodule of $A_l[l].$ It is clear
from the decomposition  (3.1) of $A_l[l]$ and Assumption II (1) and
(2) ($ii$) that every $G_F$ submodule of $A_l[l]$ is of the form
$A_l[\lambda_{i_1}] \oplus\dots\oplus  A_l[\lambda_{i_r}]$ for some
$i_1,\dots, i_r \in \{1,\dots, e\}.$ Hence if $Im\, \phi_P$ were a
proper $G_F$ submodule, we could assume that $$Im\, \phi_P \subset
A_l[\lambda_{1}] \oplus\dots\oplus  A_l[\lambda_{i-1}] \oplus
A_l[\lambda_{i+1}]\oplus\dots\oplus A_l[\lambda_{e}]$$ for some $1\leq
i \leq e.$ This implies that
$$\pi_1\dots\pi_{i-1}\pi_{i+1}\dots\pi_{e}( \sigma({1 \over
l}{\hat P}) - {1 \over l}{\hat P}) = 0 \tag{4.8}$$ for every
$\sigma\in G({\bar F}/F_l).$ The equality (4.8) takes place in
$J_{f,S_l}(T_l)$ under the (2) (i) part of Definition 3.2 (resp. in
$B({\overline F})$ under the case (2) (ii) of Definition 3.2) and it
implies that 
$$\sigma(\pi_1\dots\pi_{i-1}\pi_{i+1}\dots\pi_{e} {1
\over l}{\hat P}) = \pi_1\dots\pi_{i-1}\pi_{i+1}\dots\pi_{e}{1
\over l}{\hat P} \tag{4.9}$$ 
for every $\sigma\in G({\bar
F}/F_l).$ Hence by Proposition 2.16 (2) (resp. by Def. 3.2, of the
Mordell-Weil ${\cal O}_E$-module $\bigl\{B(L)\bigr\}$) we get
$$\pi_1\dots\pi_{i-1}\pi_{i+1}\dots\pi_{e}{1 \over l}{\hat P} \in
H^1_{f,S_{l}}(G_{F_l}; T_l)\quad (\,\in B({F_l})\,\, resp. ).
\tag{4.10}$$ 
So $\pi_1\dots\pi_{i-1}\pi_{i+1}\dots\pi_{e}{\hat P} =
0$ in the group $H^1_{f,S_{l}}(G_{F_l}; T_l)/l$ \, (in $B({F_l})/lB({F_l}),$ resp.). 
By parts (2) and (3) of the
Proposition (see also the diagram (4.6)) this implies
$\pi_1\dots\pi_{i-1}\pi_{i+1}\dots\pi_{e}{\hat P} = 0$ in the
group $B({F})/lB({F})$ in both cases. Hence there is $P_1 \in
B(F)$ such that $\pi_1\dots\pi_{i-1}\pi_{i+1}\dots\pi_{e}{\hat P}
= lP_1.$ This gives the equality
$$\pi_1\dots\pi_{i-1}\pi_{i+1}\dots\pi_{e}({\hat P} - \pi_i P_2) =
0 \tag{4.11}$$ where $P_2 = uP_1 \in B(F)$ for some $u \in {\cal
O}_E^{\times}.$ Multiplying equation (4.11) by $\pi_i$ we obtain
the equality $l(P - \pi_i P_2) = 0$ in the group $B(F).$ Since, by
assumption, $\sharp\, B(F)_{l} = 0$ we get $P = \pi_i P_2,$ hence
$P \in \lambda_i B(F)$ which contradicts the assumptions.
\qed\enddemo
\bigskip

\noindent For a given $l$ let ${\bar\rho}_{i}$ denote the
representation: $${\bar\rho}_{i} : G_{F} \rightarrow
GL(A_l[\lambda_i])$$ Similarly to the definition of $F_l$ we put
$F_i = {\bar F}^{ker {\bar\rho}_{i}}.$ In analogy with the
Definition 4.1 we introduce a homomorphism $$ \phi_i :  G_{F_i}
\rightarrow A_l[\lambda_i],$$ $$ {\phi_i}(\sigma) = {\sigma}({1
\over {\pi}_i}{\hat P}) - {1 \over {\pi}_i}{\hat P}.$$
\medskip

\noindent
\proclaim{Proposition 4.12}  We have the following properties.
\roster
\item[{1}] $H^r(G(F_i/F); A_l[{\lambda}_i])=0$ for $r\geq 0$ , all $l \in
{\cal P},$ and $1 \leq i \leq e$
except the case of trivial $G(F_i/F)$-module $A_l[\lambda_i]$
when $r = 0.$
\item[{2}] The map $H^1_{f,S_{l}}(G_F;
T_l)/{\lambda}_i \longrightarrow H^1_{f,S_{l}}(G_{F_i};
T_l)/{\lambda}_i$ is injective for all $l \in {\cal P}$ and  $1 \leq i
\leq e.$
\item[{3}] The map $B(F)/{\lambda}_i B(F)\longrightarrow
B(F_i)/{\lambda}_i B(F_i)$
is injective for all $l \in {\cal P}$  and  $1 \leq i \leq e.$
\item[{4}] Let $P\in B(F).$ If $l \in {\cal P}$ does not divide $\sharp\, B(F)_{tor}$
and $P \notin \lambda_{i}B(F),$ then the map $\phi_i$ is
surjective.
\endroster
\endproclaim
\medskip

\noindent
\demo{Proof} Proofs of (1), (2), and (3) are done in the same way as the corresponding
proofs in Proposition 4.3.
Statement (4) holds because $\phi_i$ is obviously
$G_F$ equivariant, $\phi_i$ is nontrivial
since $P \notin \lambda_{i}B(F),$
and $A_l[\lambda_i]$ is an irreducible $\Z/l[G_F]$ module
due to Assumption II.
\qed\enddemo
\medskip

\noindent
Let $P,\, Q$ be two nontorsion
elements of the group $B(F).$ Let $S_l$ be the finite set of
primes which contains primes for which $\rho_l$ is ramified and
primes over $l.$ For $v \not\in S_l$ let
$$r_v\,\,\colon\,\,
H^1_{f,S_l}(G_F;\, T_l)\rightarrow H^1(g_v;\, T_l)$$ denote the
reduction map at a prime ideal $v$ of ${\cal O}_F$.
We will
investigate the linear dependence of $P$ and $Q$ over ${\cal O}_E$ in
$B(F)$ under some local conditions for the maps $r_v,$ (see
statement of Theorem 5.1 below).
\noindent We need some additional notation. Let ${\cal P}^{\ast}$
be the set of rational primes $l \in {\cal P}$ such that $P\notin
\lambda_i B(F)$ and $Q\notin \lambda_i B(F)$ for all $1 \leq i
\leq e.$ The set ${\cal P}\setminus {\cal P}^{\ast}$ is finite,
since $B(F)$ is finitely generated ${\cal O}_E$ - module. Let ${\hat
R}\in J_{f,S_l}(T_l)$ be such that $l{\hat R}={\hat P.}$ The
element ${\hat R}$ exists by Proposition 2.14 . The
Galois group $G_{F_l}$ acts on the set $$ \{{\hat R}+t\colon\quad
t\in A_l[l]\} $$

\noindent
which is contained in $J_{f,S_l}(T_l).$  Let $N_P \subset G_{F_l}$
be the kernel of this action.
Note that $N_P$ is a normal subgroup of $G_{F_l}$ of finite index.
Define the field
$$F_l({\frac{1}{l}}{\hat P})={\bar F}^{N_P}.
$$
Let $F_l({\frac{1}{l}}{\hat Q})$ denote the corresponding field
defined for $Q.$ Observe that
$F_l({\frac{1}{l}}{\hat P})/F$ and
$F_l({\frac{1}{l}}{\hat Q})/F$ are
Galois extensions and we have isomorphisms

$$
Gal(F_l({\frac{1}{l}}{\hat P})/F)\cong H_2\rtimes
G_l\quad\quad\quad\quad
Gal(F_l({\frac{1}{l}}{\hat Q})/F)\cong H_1\rtimes G_l,
$$

\noindent
where
$$H_1=Gal(F_l({\frac{1}{l}}{\hat Q})/F_l)\quad\quad\quad\quad 
H_2=Gal(F_l({\frac{1}{l}}{\hat P})/F_l).
$$

\noindent
By Proposition 4.3 (4) the group $H_1$ (
$H_2,$ respectively) can be identified with $A_l[l]$ via the
map $\phi_Q$ ( $\phi_P,$ resp.). Put
$K = F_l({\frac{1}{l}}{\hat P})F_l({\frac{1}{l}}{\hat Q}).$
\medskip

\noindent
All fields
introduced above are displayed in the diagram below.
$$
\diagram
&K&                  \\
F_l({\frac{1}{l}}{\hat P})\urline&& F_l({\frac{1}{l}}{\hat Q})\ulline
&& \\
F({\hat R})\uline&F_l\ulline\urline&F({\hat R}^{\prime})\uline
 &&\\
&F\uline\ulline\urline&           &&
\enddiagram\tag{4.13}$$
\bigskip

\noindent
Similarly, let ${\hat R_i}\in J_{f,S_l}(T_l)$ be such that
$\pi_i {\hat R_i}={\hat P.}$ The element ${\hat R_i}$ exists by Proposition 2.14.
The Galois group $G_{F_i}$ acts on the set
$$
\{{\hat R_i}+t\colon\quad t\in A_l[\lambda_i]\}
$$

\noindent
which is contained in $J_{f,S_l}(T_l).$  Let $N_i \subset G_{F_i}$
be the kernel of this action.
Note that $N_i$ is a normal subgroup of $G_{F_i}$ of finite index.
Define the field
$$F_i({\frac{1}{\pi_i}}{\hat P})={\bar F}^{N_i}.
$$
Let $F_i({\frac{1}{\pi_i}}{\hat Q})$ denote the corresponding field
defined in the same way for $Q.$ Observe that
$F_i({\frac{1}{\pi_i}}{\hat P})/F$ and
$F_i({\frac{1}{\pi_i}}{\hat Q})/F$ are
Galois extensions and there are isomorphisms

$$
Gal(F_i({\frac{1}{\pi_i}}{\hat P})/F)\cong H_{2,i}\rtimes
G(F_i/F)\quad\quad\quad\quad
Gal(F_i({\frac{1}{\pi_i}}{\hat Q})/F)\cong H_{1,i}\rtimes G(F_i/F),
$$

\noindent
where
$$H_{1,i}=Gal(F_i({\frac{1}{\pi_i}}{\hat Q})/F_i)\quad\quad
\quad\quad H_{2,i}=Gal(F_i({\frac{1}{\pi_i}}{\hat P})/F_i).
$$

\noindent
By Proposition 4.12 (4)  the group $H_{1,i}$ (
$H_{2,i},$ respectively) can be identified with $A_l[\lambda_i]$ via the
map $\phi_i$ for $Q$ (for $P$ resp.) Put $K_i = F_i({\frac{1}{\pi_i}}{\hat
P})F_i({\frac{1}{\pi_i}}{\hat Q}).$
\medskip

\noindent
Fields
introduced above are displayed in the left diagram below.
In the right diagram  we
depicted the relevant prime ideals that will be used in the proof
of Theorem 5.1 below.
$$
\diagram
&K_i&                   && &w\quad w^{\prime}&\\
F_i({\frac{1}{\pi_i}}{\hat P})\urline&& F_i({\frac{1}{\pi_i}}{\hat Q})\ulline
&& {\beta}\urline && {\beta^{\prime}}\ulline\\
F({\hat R_i})\uline&F_i\ulline\urline&F({\hat R_i}^{\prime})\uline
 && u\uline&{}\ulline\urline&u^{\prime}\uline\\
&F\uline\ulline\urline&           &&            &v\uline\ulline\urline
\enddiagram\tag{4.14}$$
\medskip

\noindent \remark{Remark 4.15} Observe that $$F_l({\frac{1}{l}}{\hat
P}) = F_1({\frac{1}{\pi_1}}{\hat P})\dots
F_i({\frac{1}{\pi_i}}{\hat P})\dots F_e({\frac{1}{\pi_e}}{\hat
P}),$$ $$F_l({\frac{1}{l}}{\hat Q}) = F_1({\frac{1}{\pi_1}}{\hat
Q})\dots F_i({\frac{1}{\pi_i}}{\hat Q})\dots
F_e({\frac{1}{\pi_e}}{\hat Q}).$$
In addition there is an equality
$$[F({\hat R_i}): F] = [F_i({\frac{1}{\pi_i}}{\hat P}): F_i],$$
since by Proposition 4.12 (4) there are
$[F_i({\frac{1}{\pi_i}}{\hat P}): F_i]$ different imbeddings of
$F({\hat R_i})$ into ${\bar F}$ that fix $F.$ Hence from the diagram
(4.14) we find out that $F({\hat R_i}) \cap F_i = F.$
\endremark
\bigskip

\noindent
\subhead 5. The support problem for $l$-adic
representations
\endsubhead

\proclaim{Theorem 5.1}
Let ${\cal P}^{\ast}$ be the infinite set of prime
numbers introduced after the proof of Prop.4.12.
Assume that for every $l \in {\cal P}^{\ast}$ the following condition holds
in the group $H^1(g_v;\, T_l).$

\roster
\item"{}"
For every integer $m$ and for almost every $v\not\in S_l$
$$
m\,r_v({\hat P})=0\quad\quad \text{implies}\quad\quad
m\,r_v({\hat Q})=0.
$$
\endroster

\noindent
Then there exist $a \in \Bbb Z$ and $f \in {\cal O}_E$
such that  $aP + fQ = 0$ in $B(F).$
\endproclaim
\medskip

\proclaim{Lemma 5.2}
Let $H_{1,i}$ and $H_{2,i}$ be two $h$-dimensional ${\Bbb F}_l$-vector
spaces equipped with the natural action of the group
$G_i=Im\, {\bar \rho_i} \subset GL_h({\Bbb F}_l).$ Let us denote by
$\Omega_i$ the semidirect product
$(H_{1,i}\oplus H_{2,i})\rtimes G_i.$ Assume that we are given
$\sigma_i \in G_i$ such that for every $h_1\in H_{1,i}$ the element
$(h_1,\, 0,\, \sigma_i)\in (H_{1,i}\oplus \{0\})\rtimes G_i$ is
conjugate
to an element $(0,\, h_2,\, \tau_i)\in (\{0\}\oplus H_{2,i})\rtimes G_i.$
Then $1$ is not an eigenvalue of the matrix $\sigma_i.$
\endproclaim
\medskip

\noindent
\demo{Proof} cf. [C-RS], Lemma 4.2.
\qed\enddemo
\medskip

\noindent \remark{Remark 5.3} Observe, that by Assumption II,  for
every $1 \leq i\leq e$ and every $l \in {\cal P}$ there
exists a matrix $\sigma_{i}\in G_i,$ such that $1$ is an
eigenvalue of $\sigma_{i}$ with an eigenspace of dimension 1.
\endremark
\medskip

\noindent

\demo{Proof of Theorem 5.1}
We want to prove that

$$F_l({\frac{1}{l}}{\hat P})=F_l({\frac{1}{l}}{\hat Q}).\tag{5.4}$$
Hence it is enough to prove that for each $1\leq i \leq e$ we have
$$F_i({\frac{1}{\pi_i}}{\hat P}) = F_i({\frac{1}{\pi_i}}{\hat Q}).\tag{5.5}$$
Suppose this is false for some $i$. Then we observe that
$$F_i({\frac{1}{\pi_i}}{\hat P}) \cap F_i({\frac{1}{\pi_i}}{\hat Q})
= F_i,$$
since both groups $H_{1,i} = G(F_i({\frac{1}{\pi_i}}{\hat Q})/F_i)$ and
$H_{2,i} = G(F_i({\frac{1}{\pi_i}}{\hat P})/F_i)$ are irreducible
$G_i = G(F_i/F)$ modules by Assumption II (1) (ii). Hence
$$Gal(K_i/F_i) \cong H_{1,i}\oplus H_{2,i} \cong
A_l[\lambda_i] \oplus A_l[\lambda_i].\tag{5.6}$$
We need the following result.
\medskip

\noindent
\proclaim{Lemma 5.7} We have the following equality
$$K_i \cap F_l = F_i.$$
\endproclaim

\noindent
\demo{Proof} By (5.6) the group $G(K_i/F_i)$ is abelian of order
$l^{2h}.$
If $h = 1,$ then $G(F_l/F_i) \subset \prod_{j=1, j\not= i}^d
GL_1({\Bbb Z}/l)$ has order relatively prime to $l$ and it is clear
that $K_i \cap F_l = F_i.$

\noindent Now assume that $h > 1.$ We observe that $$\prod_{j=1, j\not= i}^{e}
[H_j, H_j] \subset \prod_{j=1, j\not= i}^{e} H_j \subset
G(F_l/F_i),$$ hence by Assumption II (1) (i) and (iii) the subgroup
$\prod_{j=1, j\not= i}^{e} [H_j, H_j]$ has index prime to $l$ in
$G(F_l/F_i).$ On the other hand $$\prod_{j=1, j\not= i}^{e} [H_j,
H_j] \subset [G(F_l/F_i), G(F_l/F_i)] \subset G(F_l/F_i),$$ hence
the group $G(F_l/F_i)^{ab} = G(F_l/F_i)/[G(F_l/F_i), G(F_l/F_i)]$
has order prime to $l.$ Let $K_0 = K_i \cap F_l.$ Then $K_0/F_i,$
as a subextension of $K_i/F_i,$ is abelian with order equal to
some power of $l.$ On the other hand $G(K_0/F_i)$ is a quotient of
the abelian group $G(F_l/F_i)^{ab},$ which has order prime to $l.$
This implies that the group $G(K_0/F_i)$ is trivial. Hence $K_0 =
F_i.$ \qed\enddemo
\medskip

\noindent
Let us now return to the proof of Theorem 5.1.
Consider the following tower of fields.
$$
\diagram
&K_iF_l&                  \\
K_i \urline &F_l\uline
&& \\
F_i\uline\urline &&F_1\dots F_{i-1}F_{i+1}\dots F_e\ulline
 &&\\
&F\ulline\urline&           &&
\enddiagram\tag{5.8}$$
\bigskip

\noindent We can regard $G_l = G(F_l/F)$ as the subgroup of
$\prod_{j=1}^{e} GL_h(\F_l).$ Let us pick $\sigma_{l} \in G_l$
such that $\sigma_{l}|F_i = \sigma_i$ and $\sigma_{l}|F_j =
\beta_j$ for all $j \not= i.$ Such a $\sigma_{l}$ exists by
Assumption II (1) (iv). Note that
$\sigma_l$ considered as a linear operator on the $\F_l$ vector
space $A_l[l]$ has an eigenvalue $1$ with the eigenspace of
dimension $1.$ Let $h_1\in H_{1,i}$ be an arbitrary element. Let
us pick an element of $G(K_i/F_i) \cong H_{1,i} \oplus H_{2,i}$
such that its projection onto $H_{1,i}$ is $h_1$ and its
projection onto $H_{2,i}$ is a trivial element. We denote this
element as $(h_1, 0).$ Taking into account Lemma 5.7, Remark 4.15 and
the isomorphism of Galois groups $Gal(K_i/F)\cong (H_{1,i} \oplus
H_{2,i})\rtimes G(F_i/F),$ we can define an element $\gamma \in
G(K_iF_l/F)$ such that $\gamma|K_i = (h_1, 0, \sigma_i),$
$\gamma|F({\hat R_i}) = id_{F({\hat R_i})}$  and $\gamma|F_l =
\sigma_l.$ By Chebotarev density theorem there exists a prime
${\tilde w}$ of $K_iF_l$ such that:

\roster
\item"{($i$)}"\,\, $Fr_{{\tilde w}}\, =\, \gamma\,\, \in\,\, G(K_iF_l/F),$
\item"{($ii$)}"\,\, the unique prime $v$ in $F$ below ${\tilde w}$ is not in
$S_l$ and satisfies the assumptions of Theorem 5.1.
\endroster

\noindent By the choice of prime $v$ we see that $$H^0(g_v;\,
A_l)[l] = \oplus_{j=1}^e H^0(g_v;\, A_l)[\pi_j] = H^0(g_v;\,
A_l)[\pi_i]$$ and also $H^0(g_v;\, A_l)[l] \cong \Z/l.$ Hence for
each $k \geq 1$ we have $$H^0(g_v;\, A_l)[l^k] = H^0(g_v;\,
A_l)[\pi_{i}^k]$$ which, together with finitness of $H^0(g_v;\,
A_l),$ shows that there is an $m$ such that $$H^0(g_v;\, A_l) =
H^0(g_v;\, A_l)[l^{m}] = H^0(g_v;\, A_l)[\pi_{i}^{m}]\tag{5.9}$$
and $ H^1(g_v;\, T_l)\cong H^0(g_v;\, A_l) $  is a finite, cyclic
group.

\noindent
Let $w$ ($u$ resp.) be the prime of $K_i$  ($F({\hat R_i})$ resp.) which is
over $v$ and below ${\tilde w}$ (cf. diagram (4.14). 
Consider the following commutative diagram.

$$
\CD
H^1_{f,S_l}(G_{K_i};\, T_l) @>{r_w}>> H^1(g_{w};\, T_l)\\
@AAA @AAA\\
H^1_{f,S_l}(G_{F_{i}({1 \over \pi_i}{\hat P})};\, T_l) @>{r_{\beta}}>>
H^1(g_{\beta};\, T_l)\\
@AAA @AAA\\
H^1_{f,S_l}(G_{F({\hat R_i})};\, T_l) @>{r_u}>> H^1(g_{u};\, T_l)\\
@AAA @A{\cong}AA\\
H^1_{f,S_l}(G_{F};\, T_l) @>{r_v}>> H^1(g_{v};\, T_l)
\endCD
\tag{5.10}$$

\noindent The lowest right vertical arrow in the diagram (5.10) is
an isomorphism because, by the choices we have made the prime $v$
splits in $F({\hat R_i})$ (which means that $k_v \cong k_u.$ Note that 
prime ideal
$v$ does not need to split completely in $F({\hat R_i})/F$ since this
extension is usually not Galois). The left vertical arrows
are embeddings by Proposition 2.16. Since $v$ splits in
$F({\hat R_i}),$ we have the following equality in the group
$H^1(g_{v};\, T_l)$

$$ r_v({\hat P})=\pi_i\,r_{u}({\hat R_{i}}).$$

\noindent
Let $t_v = l^{m}$
denote the order of the finite cyclic group
$H^1(g_v;\, T_l)\cong H^0(g_v;\, A_l).$  For some $c \in {\cal O}_E^{\times}$
we have
$${\frac{t_v}{l}}\,r_v({\hat
P})={\frac{t_v}{l}}\pi_i\,r_{u}({\hat R_i}) = l^{m-1}\pi_i\,r_{u}({\hat
R_i}) = c\prod_{j\not= i}\pi_j^{m-1}(\pi_i^{m}\,r_{u}({\hat
R_i})) = 0 \tag{5.11} $$
in the group $H^1(g_v;\, T_l),$ since $r_{u}({\hat R_i}) \in
H^0(g_v;\, A_l)[\pi_i^m]$ by (5.9).

\noindent
By the assumption of Theorem 5.1, equality (5.11) implies that
$$
{\frac{t_v}{l}}r_v({\hat Q})=0.\tag{5.12}
$$

\noindent Since $H^1(g_v;\, T_l)$ is cyclic, the equality (5.12)
implies that $$r_v({\hat Q}) \in l\,H^1(g_v;\, T_l).$$ This gives
$$ r_v({\hat Q})=\pi_i\,{\tilde R}^{\prime\prime}_{i}\tag{5.13} $$

\noindent 
for some ${\tilde R}^{\prime\prime}_{i} \in H^1(g_v;\,
T_l).$ By Proposition 2.14 we can find an element ${\hat
R}^{\prime\prime}_{i} \in J_{f,S_l}(T_l)$ such that $$
\pi_i\,{\hat R}^{\prime\prime}_{i}= {\hat Q}.\tag{5.14} $$ Choose
a prime $u^{\prime\prime}$ in $F({\hat R_i}^{\prime\prime})$ over
$v.$ Let $w^{\prime}$ be a prime over $u^{\prime\prime}$ in $K_i.$
Observe that, by the diagram similar to diagram 5.10 
with ${\hat P}$ and ${\hat
R_i}$ replaced by ${\hat Q}$ and
${\hat R}^{\prime\prime}_{i}$ we obtain by (5.14) that $$r_v({\hat Q})=\pi_i\,
r_{u^{\prime\prime}}({\hat R_i}^{\prime\prime}). \tag{5.15}$$

\noindent in the group $H^1(g_{u^{\prime\prime}};\, T_l),$ hence
also in  $H^1(g_{w^{\prime}};\, T_l).$  By (5.13) and (5.15) we get

$$ r_{u^{\prime\prime}}({\hat R}^{\prime\prime}_{i})- {\tilde
R}^{\prime\prime}_{i}\in A_l[\pi_i]\cap H^1(g_{u^{\prime\prime}};\,
T_l). $$

\noindent Because $A_l[\pi_i]\subset H^1_{f,S_l}(G_{K_i};\, T_l)$
(cf. proof of Lemma 2.11 and diagram (2.12)),  by Lemma 2.13 there
exists ${\hat P_0}\in H^1_{f,S_l}(G_{K_i};\, T_l)$ such that
$r_{w^{\prime}}({\hat P_0})= r_{u^{\prime\prime}}({\hat
R}^{\prime\prime}_{i})- {\tilde R}^{\prime\prime}_{i}.$ We have
the following equality

$$ r_{w^{\prime}}({\hat R}^{\prime\prime}_{i}-{\hat P_0})= {\tilde
R}^{\prime\prime}_{i}. $$ in the group $H^1(g_{w^{\prime}};\,
T_l).$

\noindent Let ${\hat R_i}^{\prime}={\hat
R}^{\prime\prime}_{i}-{\hat P_0}.$ Since $F({\hat R_i}^{\prime})
\subset F_l({\frac{1}{l}}{\hat Q})$ there is a unique prime
$u^{\prime}$ in $F({\hat R_i}^{\prime})$ below $w^{\prime}$ and
above $v.$ Of course  $r_{u^{\prime}}({\hat R_i}^{\prime})=
{\tilde R}^{\prime\prime}_{i}.$ Consider the following commutative
diagram.

$$
\CD
H^1_{f,S_l}(G_{K_i};\, T_l) @>{r_{w^{\prime}}}>> H^1(g_{w^{\prime}};\, T_l)\\
@AAA @AAA\\
H^1_{f,S_l}(G_{F_i({1 \over \pi_i}{\hat Q})};\, T_l)
@>{r_{{\beta}^{\prime}}}>> H^1(g_{{\beta}^{\prime}};\, T_l)\\
@AAA @AAA\\
H^1_{f,S_l}(G_{F({\hat R_i}^{\prime})};\, T_l) @>{r_{u^\prime}}>>
H^1(g_{u^{\prime}};\, T_l)\\
@AAA @A{=}AA\\
H^1_{f,S_l}(G_{F};\, T_l) @>{r_{v^{\prime}}}>> H^1(g_{v^{\prime}};\, T_l)
\endCD
\tag{5.16}$$

\noindent Let $Fr_{w^{\prime}} \in G(K_i/F)$ be an element of the
conjugacy class of the Frobenius element of $w^{\prime}$ over $v.$
Observe that $$Fr_{w^{\prime}}({\hat R_i}^{\prime}) = {\hat
R_i}^{\prime} + {\hat P_{0}}^{\prime}$$ for some ${\hat
P_{0}}^{\prime} \in A_l[l].$ Note that
$$Fr_{w^{\prime}}(r_{w^{\prime}}({\hat R_i}^{\prime})) =
r_{w^{\prime}}({\hat R_i}^{\prime})\tag{5.17}$$
 because
$$r_{w^{\prime}}({\hat R_i}^{\prime}) = r_{u^{\prime}}({\hat
R_i}^{\prime}) = {\tilde R}^{\prime\prime}_{i} \in H^1(g_v;
T_l).$$ On the other hand $$Fr_{w^{\prime}}(r_{w^{\prime}}({\hat
R_i}^{\prime})) = r_{w^{\prime}}(Fr_{w^{\prime}}({\hat
R_i}^{\prime})) = r_{w^{\prime}}({\hat R_i}^{\prime} + {\hat
P_{0}}^{\prime}) = r_{w^{\prime}}({\hat R}^{\prime}_{i}) +
r_{w^{\prime}}({\hat P_{0}}^{\prime}).\tag{5.18}$$ Equations (5.17)
and (5.18) show that $r_{w^{\prime}}({\hat P_{0}}^{\prime}) = 0.$
This by Lemma 2.13 implies that ${\hat P_{0}}^{\prime} = 0$ . So
$Fr_{w^{\prime}} \in G(K_i/F({\hat R_i}^{\prime})) \cong
H_{1,i}\rtimes G_i.$ Hence $Fr_{w} = (h_1,0, \sigma_i)$ is
conjugate to $Fr_{w^{\prime}} = (0, h_2, \tau_i)$ for some $h_2
\in H_{2,i}$ and $\tau_i \in G_i.$ Lemma 5.2 implies that no
eigenvalue of $\sigma_i$ is  equal to $1.$ This contradicts the
properties of $\sigma_i$ (cf. Assumption II). 
So we proved that the equality (5.5),
and consequently the equality (5.4), holds. Equality (5.4) shows
that $ker\, \phi_P = ker\, \phi_Q,$ which gives the following
commutative diagram $$ \CD 0@>>>ker({ \phi_Q})@>>> G({\bar
F}/{F_{l}})@>{\phi_Q}>> A[l] @>>> 0\\ @. @V{=}VV @V{=}VV
@V{\psi}VV\\ 0@>>>ker({ \phi_P})@>>> G({\bar
F}/{F_{l}})@>{\phi_P}>> A[l] @>>> 0\\
\endCD
\tag{5.19}$$ 
with $\psi$ a $G_l$-equivariant map. Hence due to
Assumption II (1) (v) and (2) ($ii$) (observe that (2) ($ii$)
implies that the centralizer of $G_l$ in the group $GL_d(\F_l)$ is
contained in the group of diagonal matrices $D_d \subset
GL_d(\F_l)$), it is clear, that $\psi$ as a linear operator is
represented by a block matrix of the form $$\pmatrix
b_1I_h&0&\dots&0\\
               0&b_2I_h&\dots&0\\
               \vdots&\vdots&\dots&\vdots\\
               0&0&\dots&b_eI_h
\endpmatrix
$$ for some $b_1, b_2,\dots, b_e \in \Z/l.$ Since ${\cal O}_E/(l)
\cong \prod_{j=1}^e \Z/l,$  there is a $b \in {\cal O}_E$ such
that $b$ modulo the ideal $(l{\cal O}_E)$ corresponds to the element 
$(b_1,\dots, b_e) \in
\prod_{j=1}^e \Z/l$ via this isomorphism. So diagram (5.19)
implies that $\phi_{P}=b\,\phi_Q,$ hence $\phi_{P - b Q}$ is a
trivial map. On the other hand the natural map $$ \theta\colon\,
B(F)/lB(F) \rightarrow  H^1(G_{F_l};\, T_l/l)= Hom(G_{F_l};\,A_l[l])
$$ $$\theta({X})=\phi_{X}$$ (where $\phi_{X}$ is the map from
Def. 4.1) is an injection since it can be expressed as a
composition of the injective map from Proposition 4.3 (3) and the bottom
horizontal, injective maps from diagrams: (4.5), (4.6) and (4.7).
Hence $P = b Q$ in $B(F)/lB(F).$  So the image of $P$ in
$$B_0=B(F)/\{cQ: c\in {\cal O}_E\}$$ 
is contained in the group $lB_0$ for all
primes $l \in {\cal P}^{*}$. Since by our assumption $B(F)$ and
therefore $B_0$ are finitely generated, we conclude that
$\bigcap_{l\in {\cal P}^{*}} lB_0$ is finite. Hence $a P= b Q$ for
some $a \in \Z - \{0\}$ and $b \in {\cal O}_E.$  For $f = - b$ we obtain 
$aP + f Q = 0.$ \qed\enddemo
\bigskip

\noindent
\subhead {6.  Applications}
\endsubhead

In this section we give applications
of Theorem 5.1 to the $l$-adic representations which were already 
discussed in Examples 3.3 - 3.7 .

\subsubhead  \rm {\bf  The cyclotomic character}
\endsubsubhead

\noindent
Consider the cyclotomic character
$$\rho_{l}\,\,:\,\, G({\bar F}/F) \rightarrow GL(\Z_l(1)) \cong
GL_1(\Z_{l}) \cong \Z_l^{\times},$$
(see Example 3.3).
There is a commutative diagram.
$$
\CD
{\cal O}_{F, S}^{\times}@>>> \prod_{v\not\in S_l} (k_v)^{\times}_l\\
@VVV @V=VV \\
H^1_{f,S_l}(G_F;\, \Bbb Z_l(1))@>>>
\prod_{v\not\in S_l} H^1(g_v;\, \Bbb Z_l(1))
\endCD
\tag{6.1}$$
where the left vertical arrow factors as:
$${\cal O}_{F, S}^{\times} \rightarrow {\cal O}_{F, S}^{\times}\otimes_{\Z}
\Z_l \rightarrow H^1_{f,S_l}(G_F;\, \Bbb Z_l(1)).$$
This map has finite kernel with order prime to $l.$ 
Diagram (6.1) and Theorem 5.1 applied to $\rho_l$
imply the following corollary.
\medskip

\proclaim{Corollary 6.2} Let $P,\, Q$ be two nontorsion elements of
the group ${\cal O}_{F,S}^{\times}.$ Assume that for almost every $v$ and
every integer $m$ the following condition holds $$
m\,r_v(P)=0\quad \text{in}\quad (k_v)^{\times}\quad\quad \text{implies}
\quad\quad m\,r_v(Q)=0 \quad\text{in}\quad (k_v)^{\times}. $$

\noindent
Then there exist $a, f \in \Z - \{0\}$ such that $P^a = Q^f$ in
${\cal O}_{F,S}^{\times}.$
\endproclaim

\noindent
\subsubhead \rm {\bf  K-theory of number fields}
\endsubsubhead

\noindent
Let $n$ be a positive integer. Consider the one dimensional
representation $$
\rho_l\,\,:\,\, G_F\rightarrow GL(\Z_l(n+1))\cong \Z_l^{\times}$$
which is given by the $(n+1)$-th tensor power of the cyclotomic
character. We use the notation of Example 3.4. We have the following
commutative diagram.

$$
\CD
K_{2n+1}(F)/{C_F}@>>> \prod_{v\not\in S_l}
K_{2n+1}(k_v)_l\\
@VV{\psi_{J,l}}V @V=VV \\
H^1(G_F;\, \Bbb Z_l(n+1))@>>>
\prod_{v\not\in S_l} H^1(g_v;\, \Bbb Z_l(n+1))
\endCD
\tag{6.3}$$

\noindent
Note that in this case
$$H^1(G_F;\, \Z_l(n+1))\cong H^1(G_{F,S_l};\, \Bbb Z_l(n+1))  \cong 
H^1_{f,S_l}(G_F;\, \Bbb Z_l(n+1))$$
and
$$K_{2n+1}(k_v)_l\cong H^1(g_v;\, \Bbb Z_l(n+1))\cong
H^0(g_v;\, \Q_l/\Z_l(n+1)).$$
It follows by the definition of $B(L)$ that
$$
\psi_{L,l}\quad\colon\quad B(L)\otimes\Z_l\cong H^1(G_L;\, \Bbb Z_l(n+1)).
$$
Hence as a consequence of Theorem 5.1  we get the following corollary
(cf. [BGK]):
\medskip

\proclaim{Corollary 6.4} Let $P,\, Q$ be two nontorsion elements of
the group $K_{2n+1}(F).$ Assume that for almost every $v$ and
every integer $m$ the following condition holds $$
m\,r_v(P)=0 \quad\text{in}\quad K_{2n+1}(k_v)
\quad\quad \text{implies}\quad\quad m\,r_v(Q)=0 \quad\text{in}\quad
K_{2n+1}(k_v).$$

\noindent
Then the elements $P$ and $Q$ of $K_{2n+1}(F)$ are
linearly dependent over $\Z.$
\endproclaim
\bigskip

\noindent
Theorem 5.1 and Corollary 6.4 have the following consequence
for the reduction maps
$$
r_v^{\prime}\colon\,
H_{2n+1}(K({\cal O}_F);\, \Z)\rightarrow H_{2n+1}(SL(k_v);\, \Z)
$$
defined  on the integral homology of the K-theory spectrum $K({\cal O}_F).$
\bigskip

\proclaim{Corollary 6.5}
Let $P^{\prime},\, Q^{\prime}$ be two nontorsion elements of
the group $H_{2n+1}(K({\cal O}_F);\, \Z).$ Assume that
for almost every prime ideal $v$ and for every integer $m$
the following condition holds in $ H_{2n+1}(SL(k_v);\, \Z):$
$$
m r_v^{\prime}(P^{\prime})=0
\quad\quad \text{implies}\quad\quad
m r_v^{\prime}(Q^{\prime})=0.
$$
Then the elements $P^{\prime}$ and $Q^{\prime}$ are linearly dependent in
the group $H_{2n+1}(K({\cal O}_F);\, \Z).$
\endproclaim
\medskip

\demo{Proof} Consider the following commutative diagram.
$$
\CD
K_{2n+1}({\cal O}_F) @>>> \prod_v K_{2n+1}(k_v)\\
@V{h_F}VV @V{\prod_v h_v}VV\\
H_{2n+1}(K({\cal O}_F);\, \Z) @>>> \prod_v H_{2n+1}(SL(k_v);\, \Z).
\endCD\tag{6.6}
$$

\noindent
The horizontal maps in the diagram (6.6) are induced by the
reductions at prime ideals of ${\cal O}_F.$ The vertical
maps are the Hurewicz
maps from $K$-theory to the integral homology
of the special linear group. Since the rational
Hurewicz map
$$h_F\otimes\Q\colon\quad K_{2n+1}({\cal O}_F)\otimes Q
\rightarrow H_{2n+1}(K({\cal O}_F);\, \Q)$$

\noindent is an isomorphism cf. [Bo], we can find $c,\,d\in \Z$
and nontorsion elements $P,\, Q\in K_{2n+1}({\cal O}_F),$ such
that
$$ h_F(P)=cP^{\prime}\quad\quad\text{and}\quad\quad
h_F(Q)=dQ^{\prime}. \tag{6.7} $$
Hence we can check that for every prime ideal $v$ the image
of the reduction map $r_v^{\prime}$ is contained in 
the torsion subgroup of $H_{2n+1}(SL(k_v);\, \Z).$

It follows by [A] that kernels of
the Hurewicz maps $h_F$ and $h_v,$ for any $v,$ are finite groups
of exponents which are divisible only by primes smaller than the
number ${\frac{n+1}{2}}.$ Let ${\cal P}^{\ast}$ be the set of all prime
numbers $l$ which are bigger than ${\frac{n+1}{2}}$ and
relatively prime to $cd\,\sharp C_{F}.$ Let $l \in {\cal P}^{\ast}.$
Consider the following diagram obtained from (6.6).

$$
\CD
K_{2n+1}({\cal O}_F)\otimes \Z_l @>>> \prod_v K_{2n+1}(k_v)_l\\
@V{h_F}VV @V{\prod_v h_v}VV\\
H_{2n+1}(K({\cal O}_F);\, \Z)\otimes \Z_l @>>>
\prod_{v} H_{2n+1}(SL(k_v);\, \Z)_l.
\endCD\tag{6.8}$$ 
To simplify notation we keep denoting the Hurewicz maps and the 
reduction maps in (6.8) by the same symbols as in the diagram
(6.6). Let ${\hat P}$ (${\hat Q}$ resp.) denote as before the
image of $P$ ($Q$ resp.) via the map $$K_{2n+1}({\cal O}_F)
\rightarrow (K_{2n+1}({\cal O}_F)/C_F)\otimes_{\Z} \Z_l \cong
H^1(G_F; \Z_l(n+1)).$$ Let $S_l$ denote the finite set of primes
of ${\cal O}_F$ which are over $l.$ Let $v\not\in S_l$ and assume
that $m r_v({\hat P})=0$ in the group $K_{2n+1}(k_v)_l \cong
H^1(g_v; \Z_l(n+1)).$ Since $r_v(P) = r_v({\hat P}),$ it follows
by the diagram (6.8) that $$ 0 = m h_v(r_v(P)) = m
r_v^{\prime}h_F(P)) = c\, m r_v^{\prime}(P^{\prime})$$ 
in the group $H_{2n+1}(SL(k_v); \Z)_l.$ Since $c$ is relatively prime to
$l,$ the last equality implies that $$m
r_v^{\prime}(P^{\prime})=0.$$ Since $r_v^{\prime}(P^{\prime}) \in
H_{2n+1}(SL(k_v); \Z)_{tor},$ there is a natural number $m_0$
which is prime to $l$ and such that $$m_{0} m
r_v^{\prime}(P^{\prime})=0$$ in the group $H_{2n+1}(SL(k_v); \Z).$
Hence, by assumption $$m_{0} m r_v^{\prime}(Q^{\prime})=0$$ in the
group $H_{2n+1}(SL(k_v); \Z).$ Since $m_0$ is prime to $l$ we get
$$m r_v^{\prime}(Q^{\prime})=0$$ in the group $H_{2n+1}(SL(k_v);
\Z)_l.$ We multiply the last equality by $d.$ The commutativity of
diagram (6.8) gives then the following equality in the group
$H_{2n+1}(SL(k_v); \Z)_l.$ $$ 0=m r_v^{\prime}(dQ^{\prime})=m
r_v^{\prime}(h_F(Q))= h_v(m r_v(Q))$$

\noindent Since by the choice of $l$ the map $h_v$ in the diagram
(6.8) is injective, for $v\not\in S_l,$   from the last equality
we obtain the following: $$ m r_v({\hat Q}) = m r_v(Q)=0. $$ Thus
we have checked that the elements ${\hat P}$ and ${\hat Q}$
satisfy the assumption of Theorem 5.1. Hence by Theorem 5.1, there are
$a,\, b\in \Z$ such that 
$$aP=bQ.\tag{6.9}$$ 
in the group $K_{2n+1}({\cal O}_F).$

\noindent Applying $h_F$ to equality (6.9) and using (6.7) we get
$$acP^{\prime}=bdQ^{\prime}. \qed$$
\enddemo

\subsubhead \rm {\bf   Abelian varieties}
\endsubsubhead

\noindent
Let $A/F$ be a simple abelian
variety of dimension $g$ defined over the number field $F.$
As usual $T_l = T_l(A)$ denotes the Tate module of $A.$
Consider the $l$-adic representation $$\rho_l\,\,:\,\, G_F
\rightarrow GL(T_l(A)).$$

\noindent
 We follow the notation introduced in Examples 3.5 - 3.7. 
 For any abelian variety $A/F$ there is the following 
commutative diagram
$$
\CD
A(F) @>>> \prod_{v\notin S_l} A_v(k_v)_l\\
@V{\psi_{L,l}}VV @VVV\\
H^1_{f,S_l}(G_F;\, T_l(A)) @>>> \prod_v H^1(g_v;\, T_l(A)).
\endCD\tag{6.10}
$$
$A_v$ denotes the reduction of $A$ mod $v.$ Observe that the right vertical
arrow is an injection.
Theorem 5.1, Examples 3.5, 3.6 and 3.7,  and the diagram (6.10) imply
the following corollary.
\medskip

\proclaim{Corollary 6.11}\newline Let $A$ be an abelian variety
of dimension $g \geq 1$, defined over the number field $F$ and
such that $A$ satisfies one of the following conditions:

\roster
\item"{(1)}" $A$ has the nondegenerate CM type with $End_{F}(A)\otimes \Q$
equal to a CM field $E$ such that $E^{H}\subset F $ (cf. example 3.5)
\item"{(2)}"  $A$ is a simple, principally polarized with real multiplication by
 a totally real
field $E=End_{F}(A)\otimes \Q$ such that  $E^{H}\subset F,$ 
and the field $F$ is sufficiently large
\footnote{It means that $G_l^{alg}$ is connected and 
$\bar\rho_{l}(G_F)\subset G(l)^{alg}(\F_l)$ for almost all $l.$ 
For details see the beginning of section 3 of [BGK1].}. 
We also assume that $dim\, A = he,$ where
$e=[E:\Q]$ and $h$ is odd (cf. example 3.6) or
$A$ is simple, principally polarised such that
$End_{F}(A) = \Z$ and
$dim\, A$ is equal to $2$\, or \, $6$ (cf. example 3.7 (b)).
\endroster

\noindent
Let $P,\, Q$ be two nontorsion elements of
the group $A(F).$ Assume that for almost every prime 
\, $v$ 
of \, ${\cal O}_{F}$ and
for every integer $m$ the following condition holds in $ A_v(k_v)$
$$m\,r_v(P)=0\quad\quad
\text{implies}\quad\quad m\,r_v(Q)=0.$$

\noindent
Then there exist $a \in \Z -\{0\}$ and $f \in {\cal O}_E -\{0\}$
such that $aP +fQ = 0$ in $A(F).$
\endproclaim


\bigskip

\Refs

\widestnumber\key{AAAA}

\ref\key A
\by D. Arlettaz
\paper The Hurewicz homomorphism in algebraic K-theory
\jour J. of Pure and Applied Algebra.
\vol 71
\yr 1991
\pages 1 - 12
\endref

\ref\key Ba
\by G. Banaszak
\paper Mod $p$ logarithms $log_2 3$ and $log_3 2$ differ for
infinitely many primes
\jour Annales Mathematicae Silesianae
\vol 12
\yr 1998
\pages 141 - 148
\endref

\ref\key BGK
\by G. Banaszak, W. Gajda, P. Kraso{\' n}
\paper A support problem for $K$-theory of number fields
\jour C. R. Acad. Sci. Paris S{\' e}r. 1 Math.
\vol 331 no. 3
\yr 2000
\pages 185 - 190
\endref

\ref\key BGK1
\by G. Banaszak, W. Gajda, P. Kraso{\' n}
\paper On Galois representations for abelian varieties with complex and real multiplications
\jour preprint (2002)
\endref

\ref\key BE
\by S. Bloch, E. Esnault
\paper The coniveau filtration and non-divisibility for algebraic
cycles
\jour Mathematische Annalen
\vol 304
\yr 1996
\pages 303 - 314
\endref

\ref\key BK
\by S. Bloch, K. Kato
\paper L-functions and Tamagawa numbers of motives,\, The Grothendieck
Festschrift,
\vol I
\yr 1990
\pages 333 - 400
\endref

\ref\key Bo
\by A. Borel
\paper Stable real cohomology of arithmetic groups
\jour Ann. Scient. Ec. Norm. Sup.
\vol $4^{ieme}$ s{\' e}ries 7
\yr 1974
\pages 235-272
\endref

\ref\key CF
\by J.W.S. Cassels, A. Fr\" ohlich eds. 
\book Algebraic number theory
\publ Academic Press, London and New York
\yr 1967
\endref

\ref\key C-RS
\by C.Corralez-Rodrig{\'a}{\~ n}ez, R.Schoof
\paper Support problem and its elliptic analogue
\jour Journal of Number Theory
\vol 64
\yr 1997
\pages 276-290
\endref

\ref\key C
\by  W. Chi
\paper $l$-adic and $\lambda$-adic representations associated
to abelian varieties defined over a number field
\jour American Jour. of Math.
\yr 1992
\vol 114, No. 3
\pages 315-353
\endref

\ref\key D1
\by  P. Deligne
\paper La conjecture de Weil I
\jour Publ. Math.  IHES
\vol 43
\yr 1974
\pages  273-307
\endref

\ref\key DF
\by W. Dwyer, E. Friedlander
\paper Algebraic and \' etale K-theory
\jour Trans. Amer. Math. Soc.
\vol 292
\yr 1985
\pages  247-280
\endref

\ref\key Har
\by R. Hartshorne
\book Algebraic Geometry
\publ GTA 52, Springer-Verlag
\yr 1977
\pages
\endref

\ref\key Ja
\by U. Janssen
\paper On the $l$-adic cohomology of varieties over number fields
and its Galois cohomology.\,\, Galois groups over $\Q$
\publ Math. Sci. Res. Inst. Pub. 16, Springer, New York, Berlin
\yr 1989
\pages
\endref


\ref\key La
\by S. Lang
 \book Complex Multiplication
\publ Springer Verlag
\yr 1983
\endref

\ref\key Mi1
\by J.S. Milne
\book \' Etale cohomology
\publ Princeton University Press
\yr 1980
\endref

\ref\key Mi2
\by J.S. Milne
\paper  Abelian varieties
Arithmetic Geometry G. Cornell, J.H. Silverman (eds.)
\publ Springer-Verlag
\yr 1986
\pages 103-150
\endref

\ref\key M
\by D. Mumford
\book Abelian varieties
\publ Oxford University Press
\yr 1988
\endref



\ref\key R1
\by K. A. Ribet
\paper Galois action on division points of abelian varieties with
real multiplications
\jour American Jour. of Math.
\yr 1976
\vol 98, No. 3
\pages 751-804
\endref

\ref\key R2
\by K. A. Ribet
\paper Dividing rational points of abelian varieties of CM type
\jour Compositio math.
\yr 1976
\vol 33
\pages 69-74
\endref



\ref\key S
\by A. Schinzel
\paper {O pokazatelnych sravneniach}
\jour {Matematicheskie Zapiski}
\yr 1996
\vol 2
\pages 121-126
\endref

\ref\key Sc
\by C.Schoen
\paper Complex varieties for which the Chow group mod $n$ is not finite
\jour preprint
\yr 1997
\vol
\pages
\endref

\ref\key Se1
\by J.P. Serre
\paper R\'esum\'es des cours au Coll\`ege de France
\jour Annuaire du Coll\`ege de France
\yr 1985-1986
\pages 95-100
\endref

\ref\key Se2
\by J.P. Serre
\paper Lettre \`a Daniel Bertrand du 8/6/1984
\jour Oeuvres. Collected papers. IV.
\publ Springer-Verlag, Berlin
\yr 1985 - 1998
\pages 21 - 26
\endref

\ref\key Se3
\by J.P. Serre
\paper Lettre \`a Marie-France Vign\'eras du 10/2/1986
\jour Oeuvres. Collected papers. IV.
\publ Springer-Verlag, Berlin
\yr 1985-1998
\pages 38-55
\endref



\ref\key ST
\by J.P. Serre, J. Tate
\paper Good reduction of abelian varieties
\jour Annals of Math.
\yr 1968
\vol 68
\pages 492-517
\endref

\ref\key Sil
\by J. Silverman
\paper Abelian varieties
\book The arithmetic of elliptic curves
\publ GTM 106, Springer-Verlag
\yr 1986
\endref



\ref\key T
\by J. Tate
\paper Relation between $K_2$ and Galois cohomology
\jour Invent. Math.
\vol 36
\yr 1976
\pages 257-274
\endref



\endRefs
\enddocument